\documentclass[11pt]{article}

\usepackage{amssymb,latexsym}
\usepackage{epsfig}
\usepackage{eufrak}
\usepackage{amsmath}
\usepackage{mathrsfs}
\usepackage{color}

\newtheorem{theorem}{Theorem}
\newtheorem{lemma}{Lemma}
\newtheorem{corollary}{Corollary}

\newcommand{\be}{\begin{equation}}
\newcommand{\ee}{\end{equation}}
\newcommand{\bee}{\begin{eqnarray*}}
\newcommand{\eee}{\end{eqnarray*}}
\newcommand{\bel}{\begin{eqnarray}}
\newcommand{\eel}{\end{eqnarray}}
\newcommand{\bec}{\begin{cases}}
\newcommand{\eec}{\end{cases}}
\newcommand{\bem}{\begin{bmatrix}}
\newcommand{\eem}{\end{bmatrix}}

\newcommand{\la}{\label}
\newcommand{\li}{\left}
\newcommand{\ri}{\right}

\newcommand{\ovl}{\overline}

\newcommand{\lc}{\lceil}
\newcommand{\rc}{\rceil}
\newcommand{\lf}{\lfloor}
\newcommand{\rf}{\rfloor}

\newcommand{\ep}{\epsilon}
\newcommand{\vep}{\varepsilon}

\newcommand{\Up}{\Upsilon}

\newcommand{\de}{\delta}

\newcommand{\ga}{\gamma}

\newcommand{\se}{\theta}

\newcommand{\ze}{\zeta}
\newcommand{\al}{\alpha}
\newcommand{\ba}{\beta}

\newcommand{\vro}{\varrho}
\newcommand{\ro}{\rho}

\newcommand{\om}{\omega}
\newcommand{\Om}{\Omega}

\newcommand{\f}{\frac}

\newcommand{\cd}{\cdots}

\newcommand{\qu}{\quad}
\newcommand{\qqu}{\qquad}

\newcommand{\mscr}{\mathscr}

\newcommand{\bb}{\mathbb}

\newcommand{\wh}{\widehat}
\newcommand{\wt}{\widetilde}

\newcommand{\bs}{\boldsymbol}

\newcommand{\arl}{\leftarrow}

\newcommand{\LRA}{\Longleftrightarrow}
\newcommand{\sh}{\slash}

\newcommand{\tx}{\text}

\newcommand{\iy}{\infty}

\newcommand{\pa}{\partial}

\newcommand{\bed}{\begin{description}}
\newcommand{\eed}{\end{description}}
\newcommand{\bei}{\begin{itemize}}
\newcommand{\eei}{\end{itemize}}
\newcommand{\ben}{\begin{enumerate}}
\newcommand{\een}{\end{enumerate}}
\newcommand{\bib}{\bibitem}
\newcommand{\beL}{\begin{lemma}}
\newcommand{\eeL}{\end{lemma}}
\newcommand{\beT}{\begin{theorem}}
\newcommand{\eeT}{\end{theorem}}
\newcommand{\sect}{\section}

\newcommand{\bpf}{\begin{pf}}
\newcommand{\epf}{\end{pf}}
\newcommand{\bsk}{\bigskip}
\newcommand{\bi}{\binom}

\newcommand{\pfbox}{\hfill\mbox{$\Box$}}

\newenvironment{pf}{\paragraph*{Proof{\rm.}}}{\pfbox\bigskip}

\begin{document}

\title{{\bf Inverse Sampling for Nonasymptotic Sequential Estimation of Bounded Variable Means}
\thanks{The author is currently with Department of Electrical Engineering,
Louisiana State University at Baton Rouge, LA 70803, USA, and Department of Electrical Engineering, Southern
University and A\&M College, Baton Rouge, LA 70813, USA; Email: chenxinjia@gmail.com}}

\author{Xinjia Chen}

\date{Revised on December 2, 2007}

\maketitle

\begin{abstract}

In this paper, we consider the nonasymptotic sequential estimation of means of random variables bounded in
between zero and one.  We have rigorously demonstrated that, in order to guarantee prescribed relative precision
and confidence level, it suffices to continue sampling until the sample sum is no less than a certain bound and
then take the average of samples as an estimate for the mean of the bounded random variable.  We have developed
an explicit formula and a bisection search method for the determination of such bound of sample sum, without any
knowledge of the bounded variable. Moreover, we have derived bounds for the distribution of sample size. In the
special case of Bernoulli random variables, we have established analytical and numerical methods to further
reduce the bound of sample sum and thus improve the efficiency of sampling.  Furthermore, the fallacy of
existing results are detected and analyzed.

\end{abstract}

\sect{Introduction}

In various fields of sciences and engineering, it is a frequent problem to estimate the means of bounded random
variables.  Specially, Bernoulli random variables constitute an extremely important class of bounded variables,
since the universal problem of estimating the probability of an event can be formulated as the estimation of the
mean of a Bernoulli  variable.  For examples, the problems of estimating network reliability \cite{Fishman}, the
probability of acceptable performance of uncertain systems \cite{KT} \cite{SR} and approximating probabilistic
inference in Bayesian network \cite{DL} can be cast into the framework of estimating the means of Bernoulli
variables.

Clearly, Bernoulli variables can be considered as a special class of random variable bounded in $[0,1]$.  In
many applications, one needs to estimate a quantity $\mu$ which can be bounded in $[0,1]$ after proper
operations of scaling and translation.  A typical approach is to design an experiment that produces a random
variable $X$ distributed in $[0, 1]$ with expectation $\mu$, run the experiment independently a number of times,
and use the average of the outcomes as the estimate \cite{Dagum}. This technique, referred to as {\it Monte
Carlo} method, has been applied to tackle a wide range of difficult problems. For instances, estimating
multidimensional integration, volume and counts \cite{Fishman} \cite{epde}, finding approximate solution to
enumeration problems \cite{KLM}, approximating the permanent of 0-1 valued matrices \cite{KKL}, solving the
Ising model of statistical mechanics \cite{JS2}, evaluating the bit error rate of communication systems
\cite{MH}.

Since the estimator of the mean of $X$ is obtained from finite samples of $X$ and is thus of random nature, for
the estimator to be useful, it is necessary to ensure with a sufficiently high confidence that the estimation
error is within certain margin.  The well known  Chernoff-Hoeffding bound \cite{Chernoff} \cite{Hoeffding}
asserts that if the sample size is fixed and is greater than $\f{\ln \f{2}{\de} } {2 \ep^2 }$, then, with
probability at least $1 - \de$, the sample mean approximates $\mu$ with absolute error $\ep$. Often, however,
$\mu$ is small and a good absolute error estimate of $\mu$ is typically a poor relative error approximation of
$\mu$ \cite{Dagum}. Therefore, we seek an $(\vep, \de)$ approximation for $\mu$ in the sense that the relative
error of the estimator is within a margin of relative error $\vep$ with probability at least $1 - \de$.  Since
 the mean value $\mu$ is exactly what we want to estimate, it is usually not easy to obtain reasonably tight lower bound
 for $\mu$.  For a sampling scheme with fixed sample size, a loose lower bound of $\mu$ can lead to a very
 conservative sample size.  For the most difficult and important case that no positive lower bound of $\mu$ is available,
it is not possible to guarantee prescribed relative precision and confidence level by a sampling scheme with a
fixed sample size.  This forces us to look at sampling methods with random sample sizes.

The estimation techniques based on sampling schemes without fixed sample sizes have formed a rich branch of
modern statistics under the heading of {\it sequential estimation}.  Wald provided a brief introduction to this
area in his seminal book \cite{Wald}.  Ghosh et al. offered a comprehensive exposition in \cite{Gosh}. In
particular, Nadas proposed in \cite{Naddas} a sequential sampling scheme for estimating mean values with
relative precision.  Nadas's sequential method requires no specific information on the mean value to be
estimated. However, his sampling scheme is of asymptotic nature. The confidence requirement is guaranteed only
as the margin of relative error $\vep$ tends to $0$, which implies that the actual sample size has to be
infinity. This drawback severely circumvents the application of his sampling scheme.  Due to the inherent
unknown statistical error, asymptotical methods have been criticized in some literatures (see, e.g.,
\cite{Fishman}, \cite{Hampel}, and the references therein). Especially, researchers in the areas of randomized
algorithms, controls and communication systems are very reluctant to use asymptotic methods for quantifying the
uncertainty of estimation for purpose of avoiding another level of uncertainty, namely, the unknown error of
inference (see, e.g., \cite{KT}, \cite{MH}, \cite{epde}, \cite{SR} and the references therein).  Nevertheless,
when nonasymptotic method is not available or too conservative, one has to resort to asymptotic methods.

In recent years, aimed at making Monte Carlo estimation a more efficient and rigorous method,  Dagum et al. and
Cheng have attempted to develop nonasymptotic sequential methods for estimating means of random variables
bounded in $[0,1]$. To guarantee prescribed relative precision and confidence level, Dagum et al. proposed in
\cite{Dagum} that one should continue sampling until the sample sum is no less than a threshold value.
Obviously, this is simply a generalization of the classical {\it inverse binomial sampling} \cite{H} \cite{H2}.
However, the determination of the threshold of sample sum is not trivial.  Dagum et al. provided  an explicit
formula for computing such threshold value for ensuring prescribed relative precision and confidence level.  In
\cite{Cheng}, Cheng attempted to improve the efficiency by using a smaller threshold value.

In this paper, we revisit  the sequential estimation of means of random variables bounded in $[0,1]$. We
discovered that Dagum et al. and Cheng have left major flaws in the determination of threshold of sample sum.
Specifically, the proof of Dagum et al. for their claim on the reliability of estimator is incomplete and the
gap cannot be filled by using their arguments. The proof of Cheng for his claim on the reliability of estimator
is basically incorrect.  Most importantly, we have developed a new approach to determine the smallest value of
threshold and thus make the sampling much more efficient.  An explicit formula for the threshold of sample sum
is also derived, which is substantially smaller than that of Dagum et al.  A direct consequence of our explicit
formula is that Dagum's claim can be proved as a special result of ours.  Moreover, we have derived general
bounds on the distribution of sample sizes.  Our method applies to arbitrary random variables bounded in
$[0,1]$. In the special case of Bernoulli random variables, we have developed a method to further reduce the
threshold value and thus improve the efficiency of sampling. In particular, a computational method is
established for computing the minimum threshold value when knowledge of the Bernoulli parameter is available.

The remainder of the paper is organized as follows.  In Section 2, our general theory of inverse sampling is
presented. We discuss inverse binomial sampling in Section 3. In Section 4, we illustrate an application example
in the performance evaluation of communication systems. Section 5 is the conclusion. All proofs are given in the
Appendices. The mistakes of existing works are examined in Appendices D and E.

Throughout this paper, we shall use the following notations. The expectation of a random variable is denoted by
$\bb{E}[.]$. The set of integers is denoted by $\bb{Z}$. The ceiling function and floor function are denoted
respectively by $\lc . \rc$ and $\lf . \rf$ (i.e., $\lc x \rc$ represents the smallest integer no less than $x$;
$\lf x \rf$ represents the largest integer no greater than $x$).  The left limit as $t$ tends to $0$ is denoted
as $\lim_{t \downarrow 0}$. The notation ``$\LRA$'' means ``if and only if''. The other notations will be made
clear as we proceed.

\section{General Inverse Sampling}

Let $X$ be a bounded random variable defined in a probability space $(\Om, \mscr{F}, \Pr)$ such that $0 \leq X
\leq 1$ and $\bb{E}[X] = \mu \in (0, 1)$.  We wish to estimate the mean of $X$ by using a sequence of i.i.d.
random samples $X_1, \; X_2, \; \cd$ of $X$ based on the following {\it inverse sampling} scheme:

{\it Continue sampling until the sample size reach a number $\bs{n}$ such that the sample sum $\sum_{i =
1}^{\bs{n}} X_i$ is no less than a positive number $\ga$}.

\bsk

We call this an {\it inverse sampling} scheme, since it reduces to the classical {\it inverse binomial sampling}
scheme \cite{H} \cite{H2} in the special case that $X$ is a Bernoulli random variable.

We shall consider the following two estimators for $\mu$:
\[
\wt{\bs{\mu}} = \f{\ga} { \bs{n} }, \qqu \qqu \wh{\bs{\mu}} = \f{\ga - 1} { \bs{n} - 1 }.
\]
Specially, when $X$ is a  Bernoulli random variable and $\ga$ is an integer, $\wt{\bs{\mu}}$ and $\wh{\bs{\mu}}$
are, respectively,  the {\it maximum likelihood estimator} and the {\it minimum variance unbiased estimator} for
the binomial parameter \cite{DeGroot}  \cite{H} \cite{H2}. It should be noted that $\wt{\bs{\mu}}$ is not an
unbiased estimator of the binomial parameter; the bias may be considerable for small values of $\ga$.

To control the uncertainty of estimation, for a margin of relative error $\vep \in (0,1)$ and a confidence
coefficient $\de \in (0,1)$, it is highly desirable to determine minimum $\ga$ such that {\small $\Pr  \{ \li |
\wt{\bs{\mu}} - \mu  \ri | < \vep \mu  \} > 1 - \de$} when the estimator $\wt{\bs{\mu}}$ is used, and {\small
$\Pr \{  \li |  \wh{\bs{\mu}} - \mu  \ri | < \vep \mu \} > 1 - \de$}  when the estimator $\wh{\bs{\mu}}$ is
used.

For this purpose, we have

 \beT \la{main}
Let $\vep \in (0,1)$ and $\ga > 1$. Let $X_1, \; X_2, \; \cd$ be a sequence of i.i.d. random variables defined
in a probability space $(\Om, \mscr{F}, \Pr)$ such that $0 \leq X_i \leq 1$ and $\bb{E}[X_i] = \mu \in (0, 1)$
for any positive integer $i$. Define $\wt{\bs{\mu}} = \f{\ga}{ \bs{n} }$ and $\wh{\bs{\mu}} = \f{\ga - 1}{
\bs{n} - 1 }$, where $\bs{n}$ is a random variable such that {\small $\bs{n} (\om) = \min \li \{ n \in \bb{Z}:
 \sum_{i =1}^n X_i(\om) \geq \ga \ri \}$} for any $\om \in \Om$. Define
\[
\wt{\mscr{Q}} (\vep, \ga) = (1 + \vep)^{- \ga} \exp \li ( \f{ \vep \ga  }{ 1 + \vep }  \ri ) +
 \li (  \f{ \ga - 1 + \vep } { \ga - \vep \ga } \ri )^\ga
 \exp \li ( \f{ 1 - \vep \ga - \vep }{ 1 - \vep }  \ri )
\]
and
\[
\wh{\mscr{Q}} (\vep, \ga) = (1 + \vep)^{- \ga} \exp \li ( \f{ \vep \ga  }{ 1 + \vep }  \ri ) +
 \li (  \f{ \ga - 1 } { \ga - \vep \ga   } \ri )^\ga \exp \li ( \f{ 1 - \vep \ga  }{ 1 - \vep }  \ri ).
\]
Then, the following statements hold true.

(I) $\Pr \li \{  \li |  \f{ \wt{\bs{\mu}} - \mu } { \mu } \ri | \geq \vep \ri \}  \leq  \wt{\mscr{Q}} (\vep,
\ga)$ provided that $\ga > \f{1 - \vep}{\vep}$.

(II) $\Pr \li \{  \li |  \f{ \wh{\bs{\mu}} - \mu } { \mu } \ri | \geq \vep \ri \}  \leq  \wh{\mscr{Q}} (\vep,
\ga)$ provided that $\ga > \f{1}{\vep}$.

(III) $\wt{\mscr{Q}} (\vep, \ga)$ is monotone decreasing with respect to $\ga > \f{1 - \vep}{\vep}$. Moreover,
for any $\de \in (0,1)$, there exists a unique number $\wt{\ga} > \f{1 - \vep}{\vep}$ such that $\wt{\mscr{Q}}
(\vep, \wt{\ga} ) = \de$.

(IV) $\wh{\mscr{Q}} (\vep, \ga)$ is monotone decreasing with respect to $\ga > \f{1}{\vep}$. Moreover, for any
$\de \in (0,1)$, there exists a unique number $\wh{\ga} > \f{1}{\vep}$ such that $\wh{\mscr{Q}} (\vep, \wh{\ga}
) = \de$.

(V) $(1 - \vep) \wh{\ga} < \wt{\ga} < \wh{\ga} < \f{ (1 + \vep) \ln \f{2}{\de}  } { (1 + \vep) \ln (1 + \vep) -
\vep } < \f{ (1 + \vep) \ln \f{2}{\de}  } { (2 \ln 2 - 1) \vep^2 } <  \f{ 4 (e - 2) (1 + \vep) \ln \f{2}{\de} }{
\vep^2 } $.  Moreover, {\small \be \la{tight0} \lim_{\de \to 0} \f{ \wt{\ga}  } {  \li [ \ln (1 + \vep) - \f{
\vep } { 1 + \vep } \ri ]^{-1} \ln \f{2}{\de} } = \lim_{\vep \to 0} \f{ \wt{\ga}  } {  \li [ \ln (1 + \vep) -
\f{ \vep } { 1 + \vep } \ri ]^{-1} \ln \f{2}{\de} } = 1. \ee}

(VI) For $\vro > \f{\mu}{\ga}$,  {\small \[ \Pr \li \{ \bs{n} \geq \f{\ga (1 + \vro) }{\mu} \ri \}  \leq
 \li ( 1 + \vro - \f{\mu}{\ga} \ri )^{ \li ( 1 + \vro - \f{\mu}{\ga} \ri ) \times
\f{\ga}{\mu} }  \times \li ( \f{1 - \mu} { 1 + \vro  - \f{\mu}{\ga} - \mu } \ri )^{ \li (1 + \vro - \f{\mu}{\ga}
- \mu \ri ) \times \f{\ga}{\mu} }. \]}

(VII) For $0 < \vro < 1 - \mu$, {\small \[
 \Pr \li \{ \bs{n}  \leq  \f{\ga (1 - \vro)}{\mu}  \ri \}  \leq
\li ( 1 - \vro \ri )^{ (1 - \vro) \times  \f{\ga}{\mu} }  \times  \li ( \f{1 - \mu}{1 - \vro  - \mu } \ri )^{ (1
- \vro - \mu)  \times  \f{\ga}{\mu}  }.
\]}
 \eeT

See Appendix A for a proof.  From Theorem \ref{main}, we can see that $\wt{\ga}$ and $\wh{\ga}$ can be readily
computed by a bisection search method by making use of the monotone properties and the bounds provided in (III),
(IV) and (V).

As an immediate application of Theorem \ref{main}, we can easily determine the bound (i.e., threshold value) of
sample sum without a lower bound of $\mu$.  Specially, we have

\begin{corollary}
\la{exp}
 Let $\vep, \de \in (0,1)$ and $\ga > 1$. Let $X_1, \; X_2, \; \cd$ be a sequence of i.i.d. random
variables defined in a probability space $(\Om, \mscr{F}, \Pr)$ such that $0 \leq X_i \leq 1$ and $\bb{E}[X_i] =
\mu \in (0, 1)$ for any positive integer $i$.  Define $\wt{\bs{\mu}} = \f{\ga}{ \bs{n} }$ and $\wh{\bs{\mu}} =
\f{\ga - 1}{ \bs{n} - 1 }$, where $\bs{n}$ is a random variable such that {\small $\bs{n} (\om) = \min \li \{ n
\in \bb{Z}: \sum_{i =1}^n X_i(\om) \geq \ga \ri \}$} for any $\om \in \Om$. Then, \be \la{requ}
 \Pr \li \{  \li |
\f{ \wt{\bs{\mu}} - \mu } { \mu } \ri | < \vep \ri \} > 1 - \de, \qqu \Pr \li \{  \li | \f{ \wh{\bs{\mu}} - \mu
} { \mu } \ri | < \vep \ri \} > 1 - \de \ee provided that \be \la{formu}
 \boxed{\ga > \f{ (1 + \vep) \ln \f{2}{\de}  } { (1 + \vep) \ln (1 + \vep) - \vep
} } \ee
\end{corollary}

Corollary \ref{exp} provides an {\it explicit bound} of sample sum in the inverse sampling scheme to ensure the
reliability requirements (\ref{requ}).  Actually, as can be seen from Theorem \ref{main}, an {\it implicit
bound}, $\wh{\ga}$,  makes the sample scheme more efficient while guaranteeing (\ref{requ}).  When $\vep$ or
$\de$ is small, the explicit bound is close to the implicit bound, as indicated by (\ref{tight0}).

In \cite{Dagum}, Dagum et al. claimed that, in order to ensure \be \la{dif}
 \Pr \li \{  \li | \f{ \wt{\bs{\mu}}
- \mu } { \mu } \ri | \leq \vep \ri \} \geq 1 - \de, \ee it suffices to have $\ga$ greater than $\Upsilon_1 = 1
+ \f{4 (e - 2) (1+ \vep)}{\vep^2} \ln \f{2}{\delta}$. For the same purpose of guaranteeing (\ref{dif}), Cheng
claimed in \cite{Cheng} that $\ga$ can be reduced as {\small $\bs{\al} = \f{ (1 + \vep) \ln \f{2}{\de_s} } { (1
+ \vep) \ln (1 + \vep) - \vep }$}, where $\de_s$ satisfies the equation {\small
\[ \la{formula}
 \li ( 1 - \f{\de_s}{2} \ri ) \li \{ (1 - \de_s) + \li [ 1 - 2 \li (
\f{\de_s}{2} \ri )^{ \f{1 + \vep}{1 + 2 \vep} } \ri ] \li ( \f{\de_s}{2} \ri )   + \li [  1 - 2 \li (
\f{\de_s}{2} \ri )^{ \f{1 + \vep}{1 + 3 \vep} } \ri ] \li ( \f{\de_s}{2} \ri )^2 \ri \} = \de. \] }  However, as
can be seen from our analysis in Appendices D and E,  their arguments in justification of the claims are
fundamentally flawed.

The chain of inequalities of statement (V) of Theorem \ref{main} show that our explicit bound (\ref{formu}) is
significantly smaller than $\Up_1$.  This indicates that the bound $\Up_1$, obtained by Dagum et al.,  indeed
suffices the need of ensuring (\ref{dif}), though their proof is not correct.

Although Cheng failed to prove his claim on the reliability of $\wt{\bs{\mu}}$, he obtained in \cite{Cheng} the
following useful bounds on the average sample size: \be \la{asn}
 \f{ \ga } { \mu } \leq \bb{E}[ \bs{n} ] < \f{ \ga } { \mu } + 1
 \ee by making use of the
observation that $X_1 + \cd + X_{\bs{n} - 1} < \ga \leq X_1 + \cd + X_{\bs{n}}$ and Wald's identity to conclude
that $\mu ( \bb{E}[ \bs{n}]  - 1 ) < \ga \leq \mu \bb{E}[ \bs{n} ]$ and thus (\ref{asn}).

From (\ref{asn}), it can be seen that the average sample size is almost proportional to $\ga$.  Hence, it is
reasonable to compare the efficiency of different inverse sampling schemes by their bounds of sample sum $\ga$.
For this purpose, we have plotted our explicit bound, implicit bound $\wh{\ga}$ and the bound, $\Up_1$,  of
Dagum et al. in Figs.  \ref{fig_Curve_A}--\ref{fig_Curve_D}.  It can be seen that the bound of sample sum of
Dagum et al. is too conservative and leads to a substantial waste of sampling effort.

\begin{figure}[htbp]
\centerline{\psfig{figure=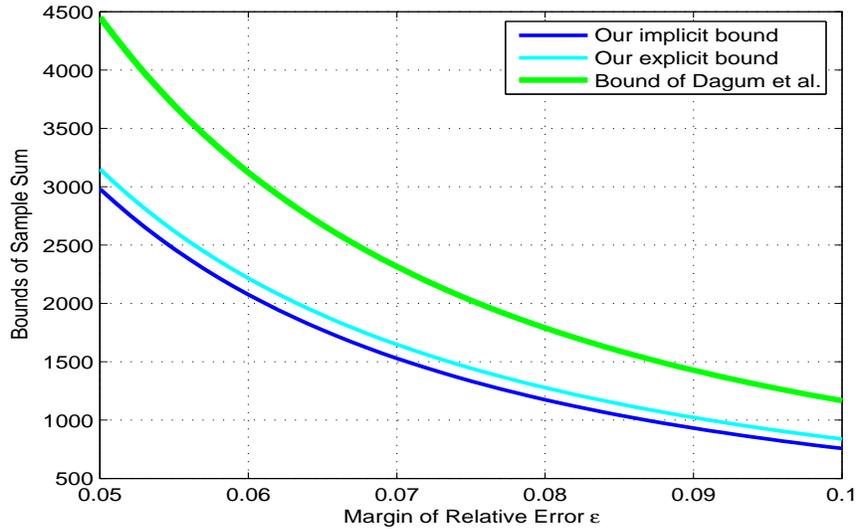, height=3in, width=5in }} \caption{Bounds of Sample Sum versus Margin of
Relative Error ($\de = 0.05$) } \la{fig_Curve_A}
\end{figure}

\begin{figure}[htbp]
\centerline{\psfig{figure=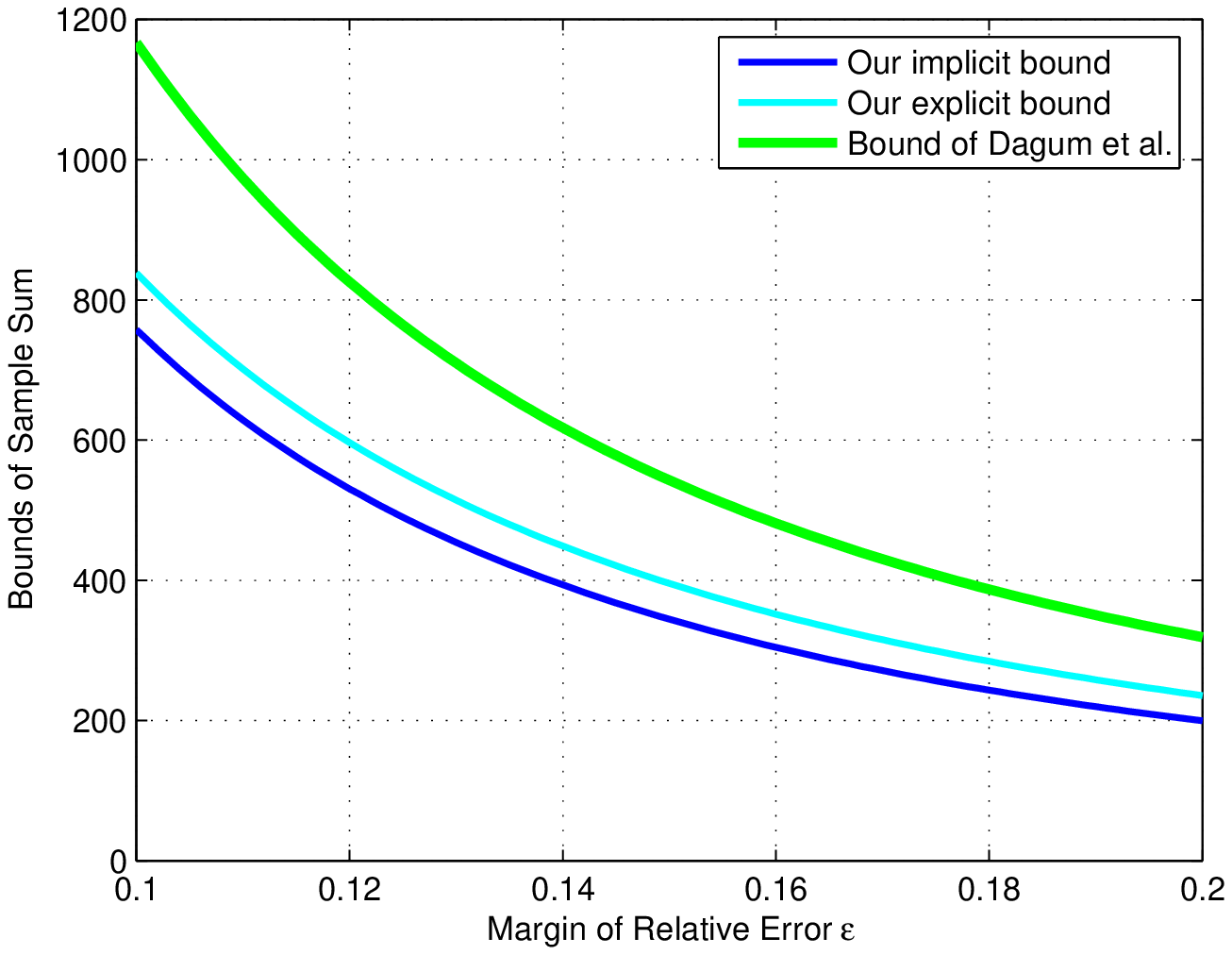, height=3in, width=5in }} \caption{Bounds of Sample Sum versus Margin of
Relative Error ($\de = 0.05$) } \la{fig_Curve_B}
\end{figure}

\begin{figure}[htbp]
\centerline{\psfig{figure=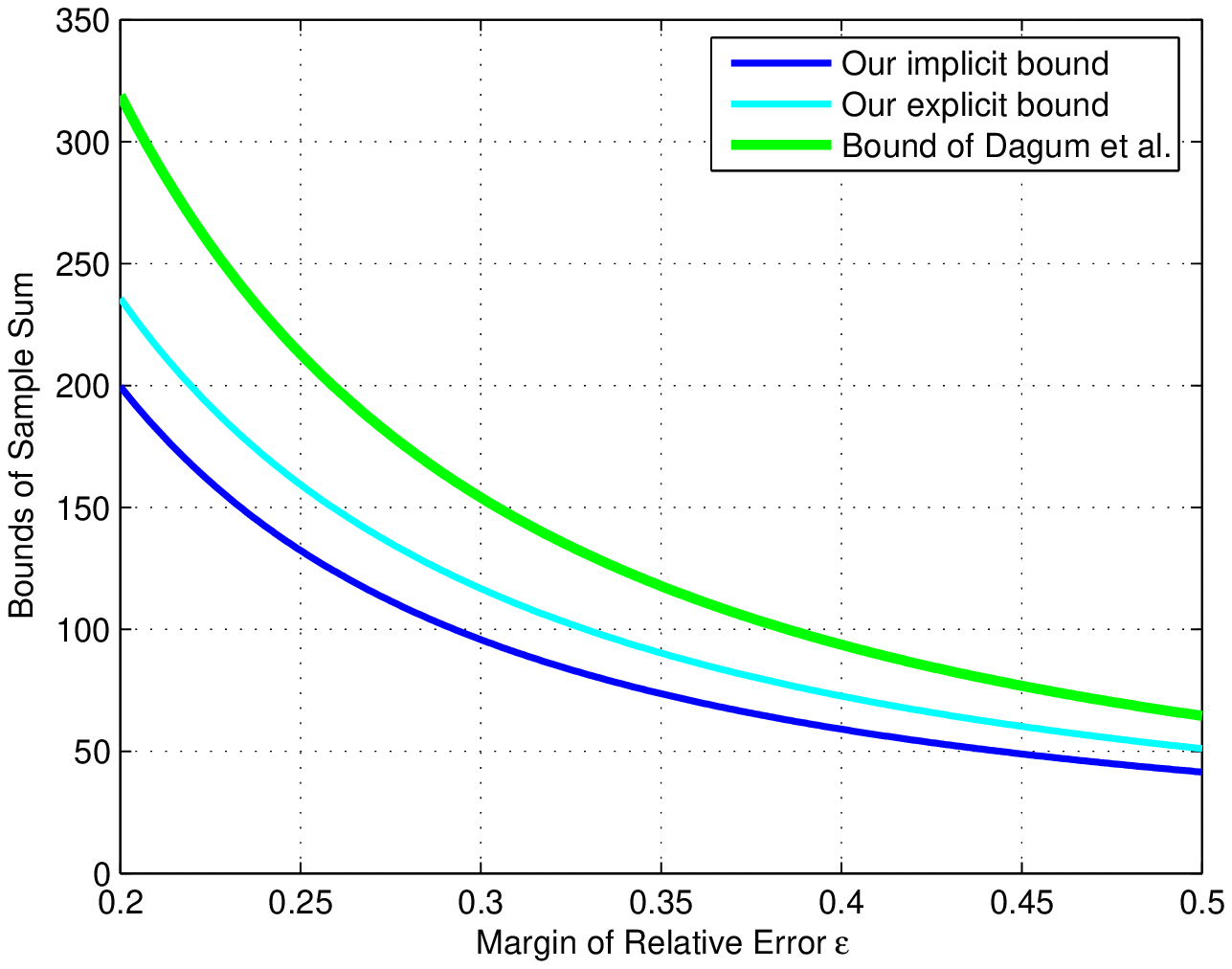, height=3in, width=5in }} \caption{Bounds of Sample Sum versus Margin of
Relative Error ($\de = 0.05$) } \la{fig_Curve_C}
\end{figure}

\begin{figure}[htbp]
\centerline{\psfig{figure=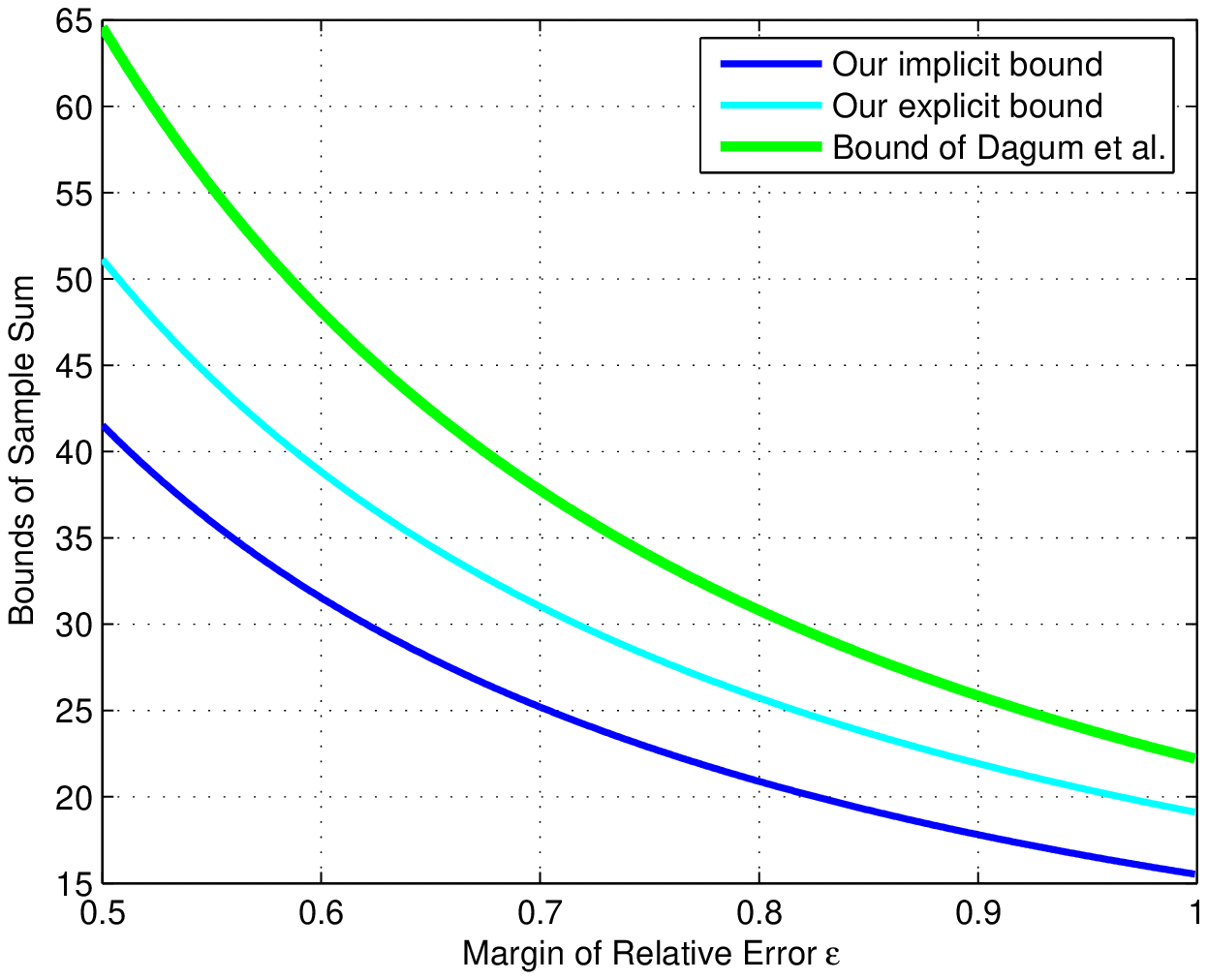, height=3in, width=5in }} \caption{Bounds of Sample Sum versus Margin of
Relative Error ($\de = 0.05$) } \la{fig_Curve_D}
\end{figure}

\section{Inverse Binomial Sampling} For the special case that $X$ is a Bernoulli random variable,
the sampling can be made more efficient.  When no knowledge of the binomial parameter is available, we have the
following results that can be used to determine the threshold value, which is smaller than its counterpart in
the general inverse sampling.

\beT \la{main_Binomial} Let $\vep \in (0,1)$. Let $X_1, \; X_2, \; \cd$ be a sequence of i.i.d. Bernoulli random
variables defined in a probability space $(\Om, \mscr{F}, \Pr)$ such that $\Pr \{ X_i = 1 \} = p \in (0,1)$ and
$\Pr \{X_i = 0 \} = 1 - p = q$  for any positive integer $i$. Define $\wt{\bs{p}}  = \f{\ga}{ \bs{n} }$, where
$\ga$ is a positive integer and $\bs{n}$ is a random variable such that {\small $\bs{n} (\om) = \min \li \{ n
\in \bb{Z}: \sum_{i =1}^n X_i(\om) = \ga \ri \}$} for any $\om \in \Om$. Then, $\Pr \li \{  \li | \f{
\wt{\bs{p}} - p } { p } \ri | \geq \vep \ri \} \leq \mscr{Q} (\vep, \ga)$ where {\small $\mscr{Q} (\vep, \ga) =
(1 + \vep)^{- \ga} \exp \li ( \f{ \vep \ga }{ 1 + \vep } \ri ) + (1 - \vep)^{- \ga} \exp \li ( - \f{ \vep \ga }{
1 - \vep }  \ri )$}, which is monotone decreasing with respect to $\ga$. Moreover, for any $\de \in (0,1)$,
there exists a unique number $\ga^*$ such that $\mscr{Q} (\vep, \ga^*) = \de$ and {\small $\max \li \{ \f{ (1 +
\vep) \ln \f{1}{\de} } { (1 + \vep) \ln (1 + \vep) - \vep }, \; \f{ (1 - \vep) \ln \f{2}{\de} } { (1 - \vep) \ln
(1 - \vep) + \vep  } \ri \} < \ga^* < \wt{\ga} < \f{ (1 + \vep) \ln \f{2}{\de} } { (1 + \vep) \ln (1 + \vep) -
\vep }$}. Furthermore, {\small
\[ \Pr \li \{ \bs{n} \geq \f{\ga (1 + \vro) }{p} \ri \} \leq \li (  \f{1 - p }{1 - p + \vro} \ri )^{ (1 - p +
\vro) \times \f{ \ga } {  p } } \li (  1 + \vro \ri )^{ (1 + \vro) \times \f{ \ga } {  p } }, \qqu \vro > 0
\]}
and {\small \[ \Pr \li \{ \bs{n} \leq \f{\ga (1 - \vro) }{p} \ri \}  \leq \li (  \f{1 - p }{1 - p - \vro} \ri
)^{ (1 - p - \vro) \times \f{ \ga } {  p } } \li (  1 - \vro \ri )^{ (1 - \vro) \times \f{ \ga } {  p } }, \qqu
0 < \vro < 1 - p.
\]}
\eeT

See Appendix B for a proof.  From Theorem \ref{main_Binomial}, it is clear that $\ga^*$ can be readily obtained
by a bisection search.

When the binomial parameter $p$ is known to be bounded in $[a, b] \subset (0,1)$, it is desirable to further
reduce the conservativeness by a computational method.  For instance, one may wish to determine the smallest
$\ga$ such that {\small $\Pr \li \{ \li | \f{ \wh{\bs{p}} - p } { p } \ri | < \vep \ri \} > 1 - \de$} for any $p
\in [a, b] \subset (0,1)$. For this purpose, an essential computational routine is to check whether a given
value of $\ga$ is large enough to ensure {\small $\Pr \li \{ \li | \f{ \wh{\bs{p}} - p } { p } \ri | < \vep \ri
\} > 1 - \de$} for any $p \in [a, b]$. At the first glance, it seems necessary to evaluate {\small $\Pr \li \{
\li | \f{ \wh{\bs{p}} - p } { p } \ri | < \vep \ri \} $} for infinite many values of $p$. Fortunately, our
following result indicates that the number of evaluations can be reduced to finite.

\beT \la{3C} Let $\vep \in (0,1)$.  Let $X_1, \; X_2, \; \cd$ be a sequence of i.i.d. Bernoulli random variables
defined in a probability space $(\Om, \mscr{F}, \Pr)$ such that $\Pr \{ X_i = 1 \} = p \in (0,1)$ and $\Pr \{X_i
= 0 \} = 1 - p = q$  for any positive integer $i$. Define $\wt{\bs{p}}  = \f{\ga}{ \bs{n} }$ and $\wh{\bs{p}}  =
\f{\ga - 1}{ \bs{n} - 1}$, where $\ga$ is a positive integer and $\bs{n}$ is a random variable such that {\small
$\bs{n} (\om) = \min \li \{ n \in \bb{Z}: \sum_{i =1}^n X_i(\om) = \ga \ri \}$} for any $\om \in \Om$. Then, the
minimum of $\Pr \li \{ \li | \f{ \wh{\bs{p}} - p } { p } \ri | < \vep \ri \}$ with respect to $p \in [a, b]
\subset (0,1)$ is achieved on the set {\small \bee &  & \{a, b \} \cup \li \{ \f{ \ga - 1 } { (1 - \vep) ( \ell
+ \ga - 1 ) } \in (a, b)
\mid \ell = 0, 1, \cd, \iy \ri \}\\
&   & \qqu \; \; \cup \li \{ \f{ \ga - 1 } { (1 + \vep) ( \ell + \ga - 1 ) } \in (a, b) \mid \ell = 0, 1, \cd,
\iy \ri \}. \eee} Similarly, the minimum of $\Pr \li \{  \li | \f{ \wt{\bs{p}} - p } { p } \ri | < \vep \ri \}$
with respect to $p \in [a, b] \subset (0,1)$ is achieved on the set {\small \bee &  & \{a, b \} \cup \li \{ \f{
\ga } { (1 - \vep) ( \ell + \ga  ) } \in (a, b) \mid \ell = 0, 1, \cd, \iy \ri \}\\
&   & \qqu \;\; \cup \li \{  \f{ \ga } { (1 + \vep) ( \ell + \ga ) } \in (a, b) \mid \ell = 0, 1, \cd, \iy \ri
\}. \eee} \eeT

\bsk

See Appendix C for a proof. The application of Theorem \ref{3C} in the computation of minimum $\ga$ is obvious.
For a fixed value of $\ga$, since the minimum of coverage probability with respect to $p \in [a, b]$ is attained
at a finite set, it can determined by a computer whether $\ga$ is large enough to ensure {\small $\Pr \li \{ \li
| \f{ \wh{\bs{p}} - p } { p } \ri | < \vep \ri \} > 1 - \de$} for any $p \in [a, b]$. Starting from $\ga = 2$,
one can find the minimum $\ga$ by gradually incrementing $\ga$ and checking whether $\ga$ is large enough.

 For $p = p_\ell = \f{ \ga - 1 } { (1 + \vep) ( \ell + \ga - 1 ) } \in (a, b)$, we have {\small \bee \Pr \li \{ \li |
\f{ \wh{\bs{p}} - p } { p } \ri | < \vep \ri \} & = & \Pr \li \{ \f{\ga - 1 }{ (1 + \vep) p } + 1 < \bs{n} <
\f{\ga - 1  }{ (1 - \vep) p }  + 1 \ri \}\\
& = & \Pr \li \{ \ell < \bs{n} - \ga < \ell + \f{ 2 \vep } { 1 -\vep } ( \ell + \ga - 1 )  \ri \}\\
& = & \sum^{\ell + \li \lc \f{ 2 \vep } { 1 - \vep } ( \ell + \ga - 1 ) \ri \rc - 1  }_{ i = \ell + 1 } \bi{ \ga
+ i - 1 } { i } p_\ell^\ga (1 - p_\ell)^i. \eee} For $p = p_\ell = \f{ \ga - 1 } { (1 - \vep) ( \ell + \ga - 1 )
} \in (a, b)$, we have {\small \bee \Pr \li \{ \li | \f{ \wh{\bs{p}} - p } { p } \ri | < \vep \ri \} & = & \Pr
\li \{ \f{\ga - 1 }{ (1 + \vep) p } + 1 <  \bs{n} < \f{\ga - 1 }{ (1 - \vep) p }  + 1 \ri
\}\\
& = & \Pr \li \{ \ell - \f{ 2 \vep } { 1 + \vep } ( \ell + \ga - 1 ) < \bs{n} - \ga < \ell  \ri \}\\
& = & \sum_{i = \max \li ( 0, \; \ell - \li \lc \f{ 2 \vep } { 1 + \vep } ( \ell + \ga - 1 ) \ri \rc + 1 \ri )
}^{ \ell - 1 } \bi{ \ga + i - 1 } { i } p_\ell^\ga (1 - p_\ell)^i. \eee}

Convenient formulas for the computation of $\Pr \li \{ \li | \f{ \wt{\bs{p}} - p } { p } \ri | < \vep \ri \}$
can be derived in a similar way.

We would like to note that the method of reducing the number of evaluations of coverage probability can also be
developed for the problems of computing minimum fixed sample sizes for the estimation of Poisson parameter,
proportions of infinite and finite populations. In this direction, we have recent research works \cite{Chen1}
\cite{Chen2} \cite{Chen3}.

\sect{An Application Example}

In this section, we shall illustrate the application of the general inverse sampling method in information
technology.   Consider the evaluation of bit error rate performance of a communication system. The stream of
bits are divided as blocks of bits with length $L > 1$. Each block is modulated as waveforms and transmitted via
a noisy channel. At the receiver side, the block of bits are recovered by demodulation.  Due to the impact of
noise, there may be incorrectly recovered bits. Let $Z$ be the number of erroneous bits. Then, $\f{Z}{L}$ is a
random variable bounded in $[0,1]$, assuming possible values $\f{\ell}{L}, \; \ell = 0, 1, \cd, L$. The bit
error rate can be defined as
\[
P_e = \bb{E} \li [ \f{Z}{L} \ri ].
\]
Since each block of bits are modulated and demodulated identically and independently, we have a sequence of
i.i.d. random variables $Z_1, \; Z_2, \; \cd$, which have the same distribution as $Z$. To estimate $P_e$, we
can continue the simulation of the modulation and demodulation process until the number of blocks reach a number
$\bs{n}$ such that
\[
\sum_{i = 1}^{\bs{n}} Z_i >  \f{ L (1 + \vep) \ln \f{2}{\de}  } { (1 + \vep) \ln (1 + \vep) - \vep }.
\]
An estimate of the bit error rate can be taken as
\[
\wh{P}_e = \f{1} { \bs{n} - 1} \li [ \f{ (1 + \vep) \ln \f{2}{\de}  } { (1 + \vep) \ln (1 + \vep) - \vep }  - 1
\ri ].
\]
Then, by our explicit formula (\ref{formu}),
\[
\Pr \li \{  \li |  \f{ \wh{P}_e - P_e } { P_e } \ri | < \vep \ri \} > 1 - \de.
\]
It should be noted that existing asymptotic estimation methods are not appropriate in this context, since the
bit error rate $P_e$ is usually very small.  Special results for Bernoulli random variables are not applicable
since $\f{Z}{L}$ is a random variable assumes $(L+ 1) > 2$ values.

\sect{Concluding Remarks}

The problem of finding relative precision estimates for means of random variables bounded in between zero and
one has numerous applications and has been studied in history for a long period of time.  Despite the lack of
rigorous justification, it was a significant progress made by a number of researchers in realizing that the
sample mean ensures the prescribed reliability once the sample sum reaches a certain threshold value.  Our main
contributions are two folds.  First, we have discovered critical mistakes exist in the determination of the
threshold value, which determines the reliability of the estimate and the efficiency of sampling. Second, we
have developed explicit formulas and computational methods to calculate the threshold value, which make the
sampling as efficient as possible, while guaranteeing prescribed relative precision and confidence level.

\appendix

\section{Proof of Theorem 1}

\beL \la{lemA1} Let $\vep \in (0,1)$ and $\ga > \f{1}{\vep}$. Let $\eta = \f{\vep \ga - 1}{\ga - 1}$ and $\ze$
be a number determined by $\f{1}{1 - \vep} = \f{1}{1 - \ze} + \f{1}{\ga}$.  Then, $0 < \eta < \ze < \vep < 1$.
\eeL

\bpf Since $0 < \vep < 1$ and $\ga > \f{1}{\vep}$, we have $\ga > 1$ and $0 < \eta = \vep - \f{1 - \vep}{\ga -
1} < \vep < 1$.

To show $0 < \ze < 1$, it suffices to exclude three possibilities.  First, $\ze \neq 1$ because $\vep \neq 1$.
 Second, $\ze
> 1$ is impossible, otherwise {\small $\f{1}{1 - \vep} < \f{1}{\ga} \LRA \ga < 1 - \vep$}, contradicting to $\ga
> \f{1}{\vep}$.  Third, $\ze \leq 0$ is impossible, otherwise {\small $\f{1}{1 - \vep} \leq 1
+ \f{1}{\ga} \LRA \ga \leq \f{ 1 - \vep } { \vep }$},  contradicting to $\ga > \f{1}{\vep}$.

Finally, since $0 < \vep < 1, \; 0 < \ze < 1$ and $\f{ 1 - \eta } {  1 - \ze } = \f{ \ga - 1 + \vep } { \ga - 1
} > 1$, we have $\ze > \eta$.

\epf

We need to use some inequalities on the function $\varphi(x) = \ln ( 1 + x ) - \f{ x }{ 1 + x }$ for $|x| < 1$.

\beL \la{lemcom} Let $\vep \in (0,1)$.  Then, $\varphi( - \vep) > \f{\vep^2}{2 (1 - \vep)} > \f{\vep^2}{2}
> \varphi(\vep) > \f{\vep^2 (2 \ln 2 - 1)} { 1 + \vep } >  0$.
 \eeL

 \bpf

To show $\varphi(\vep) < \f{ \vep^2 } { 2 }$, it suffices to note that $\varphi(\vep) = \f{ \vep^2 } { 2 } = 0$
for $\vep = 0$ and that {\small $\f{d [\varphi(\vep)  - \f{ \vep^2 } { 2 } ] } { d \vep } = \f{\vep}{(1 +
\vep)^2} - \vep < 0$} for $0 < \vep < 1$.

To show $\varphi( - \vep) > \f{\vep^2}{2 (1 - \vep)}$, it suffices to show $(1 - \vep) \ln (1 - \vep) + \vep -
\f{ \vep^2 } { 2 } > 0$ for $\vep \in (0,1)$.  This is true because the left side assumes value $0$ for $\vep =
0$ and its derivative is $ - \ln (1 - \vep) - \vep > 0$ for any $\vep \in (0,1)$.

Define $f(\vep) = (1 + \vep) \ln (1 + \vep) - \vep - (2 \ln 2 - 1) \vep^2$.
 To show {\small $\varphi(\vep) > \f{\vep^2 (2 \ln 2 - 1)} { 1 + \vep }$}, it suffices to show
$f(\vep) > 0$.  Note that $f^\prime(\vep) = \ln (1 + \vep) -  2 (2 \ln 2 - 1) \vep$ and $f^{\prime \prime}
(\vep) = \f{1}{1 + \vep} - 2 (2 \ln 2 - 1) > 0$ if $\vep < \f{ 1 }{ 2 (2 \ln 2 - 1) } - 1$.  Hence,
$f^\prime(\vep)$ is increasing for $0 < \vep < \f{ 1 }{ 2 (2 \ln 2 - 1) } - 1$ and decreasing for $\f{ 1 }{ 2 (2
\ln 2 - 1) } - 1 < \vep < 1$. Since $f^\prime(0) = 0$ and $f^\prime(1) = \ln 2 - 2 (2 \ln 2 - 1) = 2 - 3 \ln 2 <
0$, we have that there exists a unique null point $\vep^\star\in \li (  \f{ 1 }{ 2 (2 \ln 2 - 1) } - 1, 1  \ri
)$ of $f^\prime(\vep)$.  This implies that $f(\vep)$ is monotone increasing for $0 < \vep < \vep^\star$ and
monotone decreasing for $\vep^\star< \vep < 1$. Observing that $f(0) = f(1) = 0$, we can conclude that $f(\vep)
> 0$ for any $\vep \in (0,1)$. \epf

\beL \la{lemA3} Let $\vep, \; \de \in (0,1)$ and $\ga \geq \f{\ln \f{2}{\de} } { \varphi(\vep) }$.  Let $\eta =
\f{\vep \ga - 1}{ \ga - 1}$. Then, $\varphi(- \eta) > \varphi(\vep) > 0$.
 \eeL

\bpf

Since $\varphi(0) = 0$ and $\varphi^\prime (\vep)  =  \f{\vep}{(1 + \vep)^2}$,  it follows that $\varphi(\vep)
> 0$ for any $\vep \in (0,1)$.  Note that
\[
1 - \eta = 1 - \li ( \vep - \f{1 - \vep}{ \ga - 1} \ri ) = (1 - \vep) \f{\ga}{ \ga - 1}
\]
and \[ \f{\eta}{1 - \eta}  =  \f{ \vep - \f{1 - \vep}{ \ga - 1}  } { (1 - \vep) \f{\ga}{ \ga - 1} } = \f{\vep}{1
- \vep} - \f{1}{(1 - \vep) \ga}. \] Hence,
\[
\varphi ( - \eta) =  \ln \li [ (1 - \vep) \li ( 1 + \f{1}{ \ga - 1} \ri ) \ri ] + \f{\vep}{1 - \vep} - \f{1}{(1
- \vep) \ga}.
\]
By the third inequality of Lemma \ref{lemcom}, we have $\varphi(\vep) < \f{\vep^2}{2}$ for $\vep \in (0,1)$ and
thus $\ga > \f{ \ln 4 } {\vep^2}$ for any $\de \in (0,1)$.  Since $\f{1}{ \ga - 1} > 0$, using the inequality
$\ln (1 + x) \geq \f{x}{1 + x}$ with $x = \f{1}{ \ga - 1}$, we have $\ln \li ( 1 + \f{1}{ \ga - 1} \ri ) \geq
\f{ \f{1}{ \ga - 1} } { 1 + \f{1}{ \ga - 1} } = \f{1}{\ga}$ and thus
\[
\varphi ( - \eta) \geq \ln (1 - \vep) +  \f{1}{\ga} + \f{\vep}{1 - \vep} - \f{1}{(1 - \vep) \ga} = \ln (1 -
\vep) + \f{\vep}{1 - \vep} - \f{\vep}{(1 - \vep) \ga}.
\]
Define {\small $w(\vep) =  \ln (1 - \vep) + \f{\vep}{1 - \vep} - \f{\vep^3}{ (1 - \vep) \ln 4}$}.
 Applying $\ga
> \f{ \ln 4 } {\vep^2}$, we have $- \f{\vep}{(1 - \vep) \ga} >  - \f{\vep^3 }{ (1 - \vep) \ln 4}$ and
thus $\varphi ( - \eta) > w(\vep)$.  To show $\varphi (\vep) < \varphi ( - \eta)$,  it suffices to show $\varphi
(\vep) < w (\vep)$. Note that $\varphi(0) - w(0) = 0$ and {\small \bee
 \varphi^\prime (\vep) - w^\prime (\vep)
& = & \f{\vep}{(1 + \vep)^2} - \li [ \f{\vep}{(1 - \vep)^2} - \f{3 \vep^2}{ (1 - \vep) \ln 4} - \f{\vep^3}{(1 -
\vep)^2 \ln 4} \ri ]\\
&  = & \f{ \vep^2 [u(\vep) - 8 \ln 2] }{(1 + \vep)^2 (1 - \vep)^2 \ln 4 } \eee} where $u(\vep) = (1 + \vep)^2 (3
- 2 \vep)$.  Since $ u^\prime (\vep) = 2 (1 + \vep) (2 - 3 \vep) = 0$ if $\vep = \f{2}{3}$, the maximum of
$u(\vep)$ over interval $[0,1]$ must achieve at $\vep = 0, \; \f{2}{3}$ or $1$. It can be checked that $u(0) =
3, \; u(1) = 4$ and $u \li (\f{2}{3} \ri ) = \f{125}{27} < 8 \ln 2$. This shows that $u(\vep) < 8 \ln 2$, which
implies $\varphi^\prime (\vep) - w^\prime (\vep) < 0$ for any $\vep \in (0,1)$. It follows that $0 <
\varphi(\vep) < w(\vep) < \varphi ( - \eta)$.

 \epf

A classical result due to \cite{Hoeffding} is restated as Lemma \ref{Hoe} as follows.

\beL \la{Hoe} Define $\mscr{M} (z, \mu) =  \ln \li ( \f{\mu}{z} \ri ) + \li ( \f{1}{z} - 1 \ri ) \ln \li (  \f{
1 - \mu } { 1 - z } \ri )$ for $0 < z < 1$ and $0 < \mu < 1$.  Let $X_1, \cd, X_n$ be i.i.d. random variables
bounded in $[0,1]$ with common mean value $\mu \in (0,1)$. Then, {\small $\Pr \li \{ \f{\sum_{i =1}^n X_i}{n}
\geq z \ri \} \leq \exp \li ( n z \mscr{M} (z, \mu) \ri )$ } for $1
> z > \mu = \bb{E}[ X_i]$. Similarly, {\small $\Pr \li \{  \f{\sum_{i =1}^n X_i}{n} \leq z  \ri \} \leq \exp \li
(  n z \mscr{M} (z, \mu) \ri )$ } for $0 < z < \mu = \bb{E}[ X_i]$. \eeL

\beL \la{dew} $\mscr{M}( (1 + \vep) \mu, \mu)$ is monotone decreasing with respect to $\vep \in \li ( 0, \;
\f{1}{\mu} - 1 \ri )$.  Similarly, $\mscr{M}( (1 - \vep) \mu, \mu)$ is monotone decreasing with respect to $\vep
\in \li ( 0, \; 1 \ri )$.
 \eeL

 \bpf For $\vep \in \li ( 0, \;
\f{1}{\mu} - 1 \ri )$, we have $0 < (1 + \vep) \mu < 1$ and {\small \bee &   & \f{ \pa  \mscr{M}( (1 + \vep)
\mu, \mu) } { \pa \vep }\\
 & = & \f{ \pa  } { \pa \vep } \li [ \ln \li ( \f{1}{1 + \vep} \ri ) + \li (\f{1}{\mu
(1 + \vep)} - 1 \ri ) \ln \li ( \f{1 - \mu}{1 - \mu (1 + \vep) } \ri )
\ri ]\\
& = &   - \f{1}{1 + \vep} - \f{1}{\mu (1 + \vep)^2} \ln \li ( \f{1 - \mu}{1 - \mu (1 + \vep) } \ri )
 + \li (\f{1}{\mu (1 + \vep)} - 1 \ri ) \f{ \mu } { 1 - \mu (1 + \vep) }\\
 & = & - \f{1}{\mu
(1 + \vep)^2} \ln \li ( \f{1 - \mu}{1 - \mu (1 + \vep) } \ri ) < 0. \eee} Similarly, {\small $\f{ \pa  \mscr{M}(
(1 - \vep) \mu, \mu) } { \pa \vep } =  \f{1}{\mu (1 - \vep)^2} \ln \li ( \f{1 - \mu}{1 - \mu (1 - \vep) } \ri )
< 0$} for $0 < \vep < 1$. This completes the proof of the lemma. \epf

\beL \la{lembound} $\mscr{M}( (1 + \vep) \mu, \mu)$ is no greater than $- \varphi(\vep)$ for $0 < \mu < \f{1}{1
+ \vep}$. Similarly, $\mscr{M}( (1 - \vep) \mu, \mu)$ is no greater than $- \varphi(- \vep)$ for $0 < \mu < 1$.
 \eeL

\bpf

First, note that \bee &   & \lim_{\mu \to 0} \mscr{M}( (1 + \vep) \mu, \mu)\\
 & = &  \lim_{\mu \to 0} \li [ \ln
\li ( \f{1}{1 + \vep} \ri ) +
\li (\f{1}{\mu (1 + \vep)} - 1 \ri ) \ln \li ( \f{1 - \mu}{1 - \mu (1 + \vep) } \ri ) \ri ] \\
& = & \ln \li ( \f{1}{1 + \vep} \ri ) - \lim_{\mu \to 0} \ln \li ( \f{1 - \mu}{1 - \mu (1 + \vep) } \ri
) + \lim_{\mu \to 0} \f{1}{\mu (1 + \vep)} \ln \li ( \f{1 - \mu}{1 - \mu (1 + \vep) } \ri )\\
& = & \ln \li ( \f{1}{1 + \vep} \ri ) + \lim_{\mu \to 0} \f{\ln (1 - \mu) - \ln [ 1 - \mu (1 + \vep) ]
}{\mu (1 + \vep)}\\
& = &  \ln \li ( \f{1}{1 + \vep} \ri ) +  \lim_{\mu \to 0}  \f{ - \f{1}{1 - \mu} + \f{ 1 + \vep } { 1 -
\mu (1 + \vep)  }  } { 1 + \vep }\\
& = &  \ln \li ( \f{1}{1 + \vep} \ri ) +  \f{ \vep } { 1 + \vep  } = - \varphi(\vep) \eee and, similarly,
{\small $\lim_{\mu \to 0} \mscr{M}( (1 - \vep) \mu, \mu) =  \ln \li ( \f{1}{1 - \vep} \ri ) -  \f{ \vep } { 1 -
\vep } = - \varphi(- \vep)$}.  Next, we need to show that $\mscr{M}( (1 + \vep) \mu, \mu)$ is monotone
decreasing with respect to $\mu$. To
this end, note that {\small \bee &   & \f{ \pa \mscr{M}( (1 + \vep) \mu, \mu)} { \pa \mu }\\
 & = & \f{ \pa  } {
\pa \mu } \li [ \ln \li ( \f{1}{1 + \vep} \ri ) + \li (\f{1}{\mu (1 + \vep)} - 1 \ri ) \ln \li ( \f{1 - \mu}{1 -
\mu (1 + \vep) } \ri ) \ri ]\\
 & = &  \f{ \pa  } { \pa \mu } \li [ \li (\f{1}{\mu (1 + \vep)} - 1 \ri ) \ln \li ( \f{1 - \mu}{1 - \mu (1 +
\vep) } \ri ) \ri ]\\
& = &  - \f{1}{\mu^2 (1 + \vep)}  \ln \li ( \f{1 - \mu}{1 - \mu (1 + \vep) } \ri )  + \li (\f{1}{\mu (1 + \vep)}
- 1 \ri ) \li [  - \f{1}{1 - \mu} + \f{ 1 + \vep } { 1 - \mu (1 + \vep)  }  \ri ]\\
& = &  - \f{1}{\mu^2 (1 + \vep)}  \ln \li ( \f{1 - \mu}{1 - \mu (1 + \vep) } \ri ) + \li (\f{1}{\mu (1 + \vep)}
- 1 \ri ) \li (  - \f{1}{1 - \mu} \ri ) + \f{1}{\mu} \\
& = &  - \f{1}{\mu^2 (1 + \vep)}  \ln \li ( \f{1 - \mu}{1 - \mu (1 + \vep) } \ri ) - \f{1}{\mu (1 - \mu) (1 +
\vep)} + \f{1}{1 - \mu} + \f{1}{\mu}\\
& = &  - \f{1}{\mu^2 (1 + \vep)}  \ln \li ( \f{1 - \mu}{1 - \mu (1 + \vep) } \ri )  + \f{\vep}{\mu (1 - \mu) (1
+ \vep)} \leq 0 \eee} if {\small $\ln \li ( \f{1 - \mu}{1 - \mu (1 + \vep) } \ri )  \geq \f{\vep \mu }{1 -
\mu}$},  i.e., \be \la{con} \ln \li (1 - \f{\vep \mu }{1 - \mu} \ri) \leq - \f{\vep \mu }{1 - \mu}. \ee  Since
$0 < \mu < \f{1}{1 + \vep}$, we have $0 < \f{\vep \mu }{1 - \mu} < 1$.  Using the fact that $\ln (1 - x) < - x$
for any $x \in (0,1)$, we can conclude (\ref{con}) and thus establish the monotone decreasing property of
$\mscr{M}( (1 + \vep) \mu, \mu)$.

Similarly, to show that $\mscr{M}( (1 - \vep) \mu, \mu)$ is monotone decreasing with respect to $\mu$, note that
{\small \bee \f{ \pa \mscr{M}( (1 - \vep) \mu, \mu) } { \pa \mu } & = &  \f{ \pa  } { \pa \mu } \li [ \ln \li (
\f{1}{1 - \vep} \ri ) + \li (\f{1}{\mu (1 - \vep)} - 1 \ri ) \ln \li
( \f{1 - \mu}{1 - \mu (1 - \vep) } \ri ) \ri ]\\
& = &  - \f{1}{\mu^2 (1 - \vep)}  \ln \li ( \f{1 - \mu}{1 - \mu (1 - \vep) } \ri )  - \f{\vep}{\mu (1 - \mu) (1
- \vep)} \leq 0 \eee} if {\small $\ln \li ( \f{1 - \mu}{1 - \mu (1 - \vep) } \ri )  \geq - \f{\vep \mu }{1 -
\mu}$}, i.e., \be \la{conze} \ln \li (1 + \f{\vep \mu }{1 - \mu} \ri) \leq  \f{\vep \mu }{1 - \mu}. \ee
 Since $\f{\vep \mu }{1 - \mu} > 0$, using the fact that $\ln (1 + x) <  x$ for any $x \in (0,\iy)$, we can conclude
(\ref{conze}) and thus establish the monotone decreasing property of $\mscr{M}( (1 - \vep) \mu, \mu)$.

Finally, since both functions $\mscr{M}( (1 + \vep) \mu, \mu)$ and $\mscr{M}( (1 - \vep) \mu, \mu)$ are monotone
decreasing with respect to $\mu$, these two functions must be bounded from above by their corresponding limit
values as $\mu$ tends to $0$, which have been obtained at the beginning of proof.  This proves the lemma.

\epf

\beL \la{decra} {\small $\mscr{M} \li ( \f{\ga \mu}{\ga (1 + \vro) - \mu}, \mu \ri )$} is monotone decreasing
with respect to {\small $\vro > \f{ \mu } { \ga }$}. \eeL

\bpf

Let {\small $z = \f{\ga \mu}{\ga (1 + \vro)  - \mu}$} and $m = \f{\ga}{z}$.  For $\vro > \f{ \mu } { \ga }$, we
have $0 < z < \mu, \; \f{ \pa m } { \pa \vro } > 0$ and  {\small \bee \f{ \pa [\ga \mscr{M}(z, \mu)] } { \pa
\vro } & = & \ga \li [- \f{1}{z} \f{ \pa z } { \pa \vro } + \li (\f{1}{z} - 1 \ri ) \f{1}{1 - z} \f{ \pa z } {
\pa \vro } - \f{1}{z^2} \f{ \pa z } {
\pa \vro } \ln \li ( \f{1 - \mu}{1 - z} \ri ) \ri ]\\
& = & - \f{\ga}{z^2} \f{ \pa z } { \pa \vro } \ln \li ( \f{1 - \mu}{1 - z} \ri ) =  - \f{m \ga}{m z^2} \f{ \pa z } { \pa \vro } \ln \li ( \f{1 - \mu}{1 - z} \ri )\\
& = & -\f{m }{ z } \f{ \pa z } { \pa \vro } \ln \li ( \f{1 - \mu}{1 - z} \ri ) =  \f{ \pa m } { \pa \vro } \ln
\li ( \f{1 - \mu}{1 - z} \ri ) < 0, \eee} which implies that $\mscr{M}(z, \mu)$ is monotone decreasing with
respect to {\small $\vro > \f{ \mu } { \ga }$}.  This proves the lemma.

\epf

We are now in position to prove Theorem 1.

\subsection{Proof of (I)} By the definition of the estimator $\wt{\bs{\mu}} = \f{\ga}{ \bs{n} }$,
\[
\Pr \li \{  \li |  \f{ \wt{\bs{\mu}} - \mu } { \mu } \ri | \geq \vep \ri \} =  \Pr \li \{  \bs{n} \leq \f{ \ga }
{ \mu (1 + \vep) } \ri \} + \Pr \li \{  \bs{n} \geq \f{ \ga } { \mu (1 - \vep) } \ri \}.
\]
Hence, we shall derive upper bounds for the tail probabilities {\small $\Pr \li \{  \bs{n} \leq \f{ \ga } { \mu
(1 + \vep) } \ri \}$} and {\small $\Pr \li \{  \bs{n} \geq \f{ \ga } { \mu (1 - \vep) } \ri \}$}.

We first bound {\small $\Pr \li \{ \bs{n} \leq \f{ \ga } { \mu (1 + \vep) } \ri \}$}. Since $\bs{n}$ is an
integer, we have
\[
 \Pr \li \{  \bs{n} \leq \f{ \ga } { \mu (1 +
\vep) } \ri \} = \Pr \li \{  \bs{n} \leq \li \lf \f{ \ga } { \mu (1 + \vep) } \ri \rf \ri \} = \Pr \li \{ \bs{n}
\leq \f{ \ga } { \mu (1 + \vep^*) } \ri \}
\]
where $\vep^*$ is a number such that {\small $\f{ \ga } { \mu (1 + \vep^*) } = \li \lf \f{ \ga } { \mu (1 +
\vep) } \ri \rf$}.  Clearly, {\small $\vep^* = \f{ \ga } { \mu  \li \lf \f{ \ga } { \mu (1 + \vep) } \ri \rf} -
1 \geq \vep > 0$}. For simplicity of notation, let {\small $m = \f{ \ga } { \mu (1 + \vep^*) }$}. Since $m$ is a
nonnegative integer, it can be zero or a natural number.  If $m = 0$, then
\[
\Pr \li \{  \bs{n} \leq \f{ \ga } { \mu (1 + \vep) } \ri \}  = \Pr \{ \bs{n} \leq m \} = 0 < \exp ( - \ga
\varphi(\vep) ).
\]
Otherwise if $m \geq 1$, then
\[
\Pr \li \{ \bs{n} \leq \f{ \ga } { \mu (1 + \vep^*) } \ri \} = \Pr \{ \bs{n} \leq m \} = \Pr \{ X_1 + \cd + X_m
\geq \ga \} = \Pr \{ \ovl{X} \geq z \},
\]
where {\small $\ovl{X} = \f{ \sum_{i=1}^m X_i } { m }$} and {\small $z = \f{\ga}{ m } = \mu (1 + \vep^*) >
\mu$}. Now we shall consider three cases.

\bed

\item  (i): In the case of $z > 1$, we have {\small $\Pr \{ \ovl{X} \geq z \} \leq \Pr \li \{ \sum_{i = 1}^m X_i > m
\ri \} = 0 < \exp ( - \ga \varphi(\vep) )$}.

\item (ii): In the case of $z = 1$, we have $\mu = \f{1}{1 + \vep^*}, \; m = \ga$ and {\small \bee \Pr \{ \ovl{X} \geq z
\} & = & \Pr \li \{ \sum_{i = 1}^m X_i = m \ri \} = \prod_{i = 1}^m \Pr \{ X_i
= 1 \} \leq  \prod_{i = 1}^m \bb{E}[ X_i ] = \mu^m\\
& = & \li ( \f{ 1 } {1 + \vep^* }  \ri )^\ga \leq \li ( \f{ 1 } {1 + \vep}  \ri )^\ga \leq \exp ( - \ga
\varphi(\vep) ). \eee}

\item (iii): In the case of  $\mu < z < 1$, by Lemma \ref{Hoe}, we have
\[
\Pr \{ \ovl{X} \geq z \} \leq \exp( m z \mscr{M} (z, \mu) ) = \exp( \ga \mscr{M} ((1 + \vep^*) \mu, \mu) ).
\]
Since $\vep^* \geq \vep$,  it must be true that $\mu (1 + \vep) \leq \mu (1 + \vep^*) < 1$ and that $\mscr{M}
((1 + \vep^*) \mu, \mu) \leq \mscr{M} ((1 + \vep) \mu, \mu)$ as a result of Lemma \ref{dew}. Hence, \be
\la{remm} \Pr \li \{  \bs{n} \leq \f{ \ga } { \mu (1 + \vep) } \ri \} = \Pr \{ \ovl{X} \geq z \} \leq \exp( \ga
\mscr{M} ((1 + \vep) \mu, \mu) ). \ee  By Lemma \ref{lembound}, we have $\mscr{M} ((1 + \vep) \mu, \mu) \leq  -
\varphi (\vep)$. It follows that {\small
\[
\Pr \li \{  \bs{n} \leq \f{ \ga } { \mu (1 + \vep) } \ri \} = \Pr \{ \ovl{X} \geq z \} \leq
 \exp ( - \ga \varphi ( \vep ) ).
 \]
 }

\eed

 Therefore, we have shown
 \be
 \la{taila}
\Pr \li \{  \bs{n} \leq \f{ \ga } { \mu (1 + \vep) } \ri \}  \leq
 \exp ( - \ga \varphi ( \vep ) )
 \ee
 for all cases.

\bsk

 We now bound $\Pr \li \{  \bs{n} \geq \f{ \ga } { \mu (1 - \vep) } \ri \}$. Since $\bs{n}$ is an integer,
we have
\[
\Pr \li \{  \bs{n} \geq \f{ \ga } { \mu (1 - \vep) } \ri \} =
 \Pr \li \{  \bs{n} \geq \li \lc \f{ \ga } { \mu (1
-\vep) } \ri \rc \ri \} = \Pr \li \{  \bs{n} > \li \lc \f{ \ga } { \mu (1 -\vep) } - 1 \ri \rc \ri \}.
\]
Let $\ze$ be a number such that {\small $\f{1}{1 -\vep} = \f{1}{ 1 - \ze} + \f{1}{\ga}$}.  Then, $0 < \ze <
\vep$ as a result of Lemma \ref{lemA1}.  By the definition of $\ze$, we have {\small $\f{ \ga } { \mu (1 - \vep)
} - 1
> \f{\ga}{\mu} \li ( \f{1}{1 -\vep} - \f{1}{\ga} \ri ) = \f{ \ga } { \mu (1 - \ze) }$} for any $\mu \in (0,1)$. Hence,
\[
\Pr \li \{  \bs{n} > \li \lc \f{ \ga } { \mu (1 -\vep) } - 1 \ri \rc \ri \} \leq \Pr \li \{  \bs{n} > \li \lc
\f{ \ga } { \mu (1 - \ze ) }  \ri \rc \ri \} = \Pr \li \{  \bs{n} >  \f{ \ga } { \mu (1 - \ze^* ) }  \ri \}
\]
with $\ze^*$ satisfying {\small $\f{ \ga } { \mu (1 - \ze^* ) }  = \li \lc \f{ \ga } { \mu (1 - \ze ) } \ri
\rc$}. Clearly, $1 > \ze^* \geq \ze > 0$.  Let {\small $m = \f{ \ga } { \mu (1 - \ze^* ) }$}. Then, $m$ is a
positive integer and
\[
\Pr \{  \bs{n} > m \} = \Pr \{ X_1 + \cd + X_m < \ga \} = \Pr \{ \ovl{X} < z \}
\]
where {\small $\ovl{X} = \f{ \sum_{i=1}^m X_i } { m }$} and $z = (1 - \ze^*) \mu$.  Applying Lemma \ref{Hoe}, we
have {\small \bee \Pr \li \{  \bs{n} > \li \lc \f{ \ga }
 { \mu (1 - \ze) } \ri \rc  \ri \}  & = & \Pr \{ \ovl{X} < z \} \leq  \exp ( m z
\mscr{M} ( z, \mu) )\\
&  = & \exp ( \ga \mscr{M} ( (1 - \ze^*) \mu, \mu) ). \eee} Note that $\mscr{M} ( (1 - \ze^*) \mu, \mu) \leq
\mscr{M} ( (1 - \ze) \mu, \mu)$ as a result of $1 > \ze^* \geq \ze > 0$ and Lemma \ref{dew}. Hence,
\[
\Pr \{ \ovl{X} < z \} \leq \exp ( \ga \mscr{M} ( (1 - \ze) \mu, \mu) ).
\]
By Lemma \ref{lembound}, we have $\mscr{M} ((1 - \ze) \mu, \mu) \leq  - \varphi ( - \ze)$.  It follows that
{\small \be \la{tailb} \Pr \li \{  \bs{n} > \li \lc \f{ \ga }
 { \mu (1 - \ze) } \ri \rc  \ri \} \leq \exp ( - \ga \varphi ( - \ze) ). \ee } Thus, we have bounds for the two tail
probabilities as follows: {\small \[ \Pr \li \{  \bs{n} \leq \f{ \ga } { \mu (1 + \vep) } \ri \} \leq \exp ( -
\ga \varphi ( \vep) ), \qu \Pr \li \{ \bs{n} \geq \f{ \ga } { \mu (1 - \vep) } \ri \} \leq \exp ( - \ga \varphi
( - \ze) ).
\]}
It follows that \bee \Pr \li \{ \li | \f{ \wt{\bs{\mu}} - \mu } { \mu }  \ri | \geq \vep \ri \}  = \Pr \li \{
\bs{n} \leq \f{ \ga }
{ \mu (1 + \vep) } \ri \} + \Pr \li \{ \bs{n} \geq \f{ \ga } { \mu (1 - \vep) } \ri \}\\
\leq \exp ( - \ga \varphi ( \vep) ) + \exp ( - \ga \varphi ( - \ze) ) = \wt{\mscr{Q}} (\vep, \ga), \eee where we
have used the definitions of $\ze$ and $\varphi(.)$ in the last equality.     This completes the proof of
statement (I).

\subsection{Proof of (II)}

By the definition of  the estimator $\wh{\bs{\mu}}= \f{\ga - 1} { \bs{n} - 1 }$, we have \[ \Pr \li \{ \li | \f{
\wh{\bs{\mu}} - \mu } { \mu }  \ri | \geq \vep \ri \} = \Pr \li \{ \bs{n} \leq 1 + \f{ \ga - 1} { (1 + \vep) \mu
}  \ri \} + \Pr \li \{  \bs{n} \geq 1 + \f{ \ga - 1} { (1 - \vep) \mu } \ri \}.
\]
To bound {\small $\Pr \li \{ \bs{n} \leq 1 + \f{ \ga - 1} { (1 + \vep) \mu }  \ri \}$}, we shall consider two
cases.

\bed

\item (i): In the case of $(1 + \vep) \mu > 1$, we have {\small $\Pr \li \{ \bs{n} \leq 1 + \f{ \ga - 1} { (1 + \vep)
\mu } \ri \} \leq \Pr \{ \bs{n} < \ga \} = 0 \leq \Pr \li \{ \bs{n} \leq \f{ \ga } { (1 + \vep) \mu }  \ri \}$}
because $\bs{n} \geq \sum_{i = 1}^{\bs{n}} X_i \geq \ga$ is always true.

\item (ii): In the case of $(1 + \vep) \mu \leq 1$, we have {\small $\Pr \li \{ \bs{n} \leq 1 + \f{ \ga - 1} { (1 +
\vep) \mu } \ri \} \leq \Pr \li \{ \bs{n} \leq \f{ \ga } { (1 + \vep) \mu } \ri \}$}.

\eed

Therefore, in both cases, we have {\small $\Pr \li \{ \bs{n} \leq 1 + \f{ \ga - 1} { (1 + \vep) \mu } \ri \}
\leq \Pr \li \{ \bs{n} \leq \f{ \ga } { (1 + \vep) \mu } \ri \}$} and, by virtue of (\ref{taila}),
 {\small
\be \la{ineup}
 \Pr \li \{ \bs{n} \leq 1 + \f{ \ga - 1} { (1 + \vep) \mu } \ri \} \leq \exp ( - \ga \varphi
(\vep) ). \ee}

To bound {\small $\Pr \li \{  \bs{n} \geq 1 + \f{ \ga - 1} { (1 - \vep) \mu } \ri \}$}, let {\small $\eta = \f{
\vep \ga - 1  } { \ga - 1 }$} and note that {\small \bel \Pr \li \{ \bs{n} \geq 1 + \f{ \ga - 1} { (1 - \vep)
\mu } \ri \} & = & \Pr \li \{ \bs{n} \geq 1 + \li \lc \f{
\ga - 1} { (1 - \vep) \mu } \ri \rc \ri \} \nonumber\\
&  = & \Pr \li \{  \bs{n} > \li \lc \f{ \ga - 1} { (1 - \vep) \mu } \ri \rc \ri \} \nonumber\\
&  = & \Pr \li \{  \bs{n} > \li
\lc \f{ \ga } { (1 - \eta) \mu } \ri \rc \ri \} \nonumber\\
& \leq & \exp ( - \ga \varphi ( - \eta) ) \la{ineqlow} \eel} where (\ref{ineqlow}) follows from a similar method
as proving (\ref{tailb}).

Combining (\ref{ineup}), (\ref{ineqlow}) and invoking the definitions of $\eta$ and $\varphi(.)$ yields {\small
$\Pr \li \{ \li | \f{ \wh{\bs{\mu}} - \mu } { \mu } \ri | \geq \vep \ri \} \leq \exp ( - \ga \varphi (\vep) ) +
\exp ( - \ga \varphi ( - \eta) ) = \wh{\mscr{Q}} (\vep, \ga)$}. This completes the proof of statement (II).

\subsection{Proof of (III)}

Note that {\small $\wt{\mscr{Q}} (\vep, \ga) = \exp ( - \ga \varphi(\vep) ) + \exp ( - \ga \varphi(- \ze) )$},
where $\ze$ is determined by {\small $\f{1}{1 - \vep} = \f{1}{1 - \ze} + \f{1}{\ga}$}. By the chain rule of
differentiation, we have {\small \[ \f{ \pa [\wt{\mscr{Q}} (\vep, \ga) ] } { \pa \ga } = - \varphi(\vep) \exp (
- \ga \varphi(\vep) ) - \exp ( - \ga \varphi(- \ze) ) \li [ \varphi(- \ze) + \ga \f{d \varphi (-\ze)}{d \ze} \f{
\pa \ze } { \pa \ga} \ri ] < 0
\]}
 by observing that $\varphi(\vep) > 0, \; \varphi(- \ze) > 0, \; \f{d \varphi (-\ze)}{d \ze} > 0$ and $\f{ \pa \ze  } { \pa \ga} >
 0$.  This proves that $\wt{\mscr{Q}} (\vep, \ga)$ is monotone decreasing with respect to $\ga > \f{1 - \vep}{\vep}$.

The existence and uniqueness of $\wt{\ga}$ in interval {\small $\li (\f{1 - \vep}{\vep}, \iy \ri )$} can be seen
by the monotone decreasing property of $\wt{\mscr{Q}} (\vep, \ga)$ with respect to $\ga > \f{1 - \vep}{\vep}$
and the fact that {\small \[ \lim_{\ga \to \iy} \wt{\mscr{Q}} (\vep, \ga) = 0, \qu \lim_{\ga \to \f{1 -
\vep}{\vep}} \wt{\mscr{Q}} (\vep, \ga)
> \lim_{\ga \to \f{1 - \vep}{\vep}} \li (  \f{ \ga - 1 + \vep } { \ga - \vep \ga } \ri )^\ga
 \exp \li ( \f{ 1 - \vep \ga - \vep }{ 1 - \vep }  \ri ) = 1.
 \]}

\subsection{Proof of (IV)}

Note that {\small $\wh{\mscr{Q}} (\vep, \ga) = \exp ( - \ga \varphi(\vep) ) + \exp ( - \ga \varphi(- \eta) )$},
where {\small $\eta = \f{ \vep \ga - 1 } { \ga - 1 }$}.  By the chain rule of differentiation, we have {\small
\[ \f{ \pa [\wh{\mscr{Q}} (\vep, \ga) ]  } { \pa \ga  } = - \varphi(\vep) \exp ( - \ga \varphi(\vep) ) - \exp (
- \ga \varphi(- \eta) ) \li [ \varphi(- \eta)  + \ga  \f{d \varphi (-\eta)}{d \eta}  \f{ \pa \eta } { \pa \ga}
\ri ] < 0
 \]}
 by observing that $\varphi(\vep) > 0, \; \varphi(- \eta) > 0, \; \f{d \varphi (-\eta)}{d \eta}  > 0$ and $ \f{ \pa \eta  } { \pa \ga} >
 0$.  This proves that $\wh{\mscr{Q}} (\vep, \ga)$ is monotone decreasing with respect to $\ga > \f{1}{\vep}$.

The existence and uniqueness of $\wh{\ga}$ in interval {\small $\li (\f{1}{\vep}, \iy \ri )$} can be seen by the
monotone decreasing property of $\wh{\mscr{Q}} (\vep, \ga)$ with respect to $\ga > \f{1}{\vep}$ and the fact
that {\small
\[ \lim_{\ga \to \iy} \wh{\mscr{Q}} (\vep, \ga) = 0, \qqu \lim_{\ga \to \f{1}{\vep}} \wh{\mscr{Q}} (\vep, \ga) >
\lim_{\ga \to \f{1}{\vep}} \li (  \f{ \ga - 1 } { \ga - \vep \ga   } \ri )^\ga \exp \li ( \f{ 1 - \vep \ga  }{ 1
- \vep }  \ri ) = 1.
\]}

\subsection{Proof of (V)}

First, we shall show the upper bounds for $\wh{\ga}$. For this purpose, note that, for $\ga \geq \f{ \ln
\f{2}{\de} } { \varphi(\vep) }$, we have $0 < \varphi(\vep) < \varphi(- \eta)$ as a result of Lemma \ref{lemA3}.
Hence, $\wh{\mscr{Q}} (\vep, \ga) < 2 \exp ( - \ga \varphi(\vep) ) \leq \de$ for $\ga \geq \f{ \ln \f{2}{\de} }
{ \varphi(\vep) }$.  By the third inequality of Lemma \ref{lemcom}, we have $\f{ \ln \f{2}{\de} } {
\varphi(\vep) } > \f{ 2 \ln \f{2}{\de} } { \vep^2 } > \f{1}{\vep}$.  Since $\wh{\mscr{Q}} (\vep, \wh{\ga}) =
\de$ and $\wh{\mscr{Q}} (\vep, \ga)$ is monotone decreasing with respect to $\ga > \f{1}{\vep}$, it must be true
that $\wh{\ga} < \f{ \ln \f{2}{\de} } { \varphi(\vep) }$. Applying Lemma \ref{lemcom}, we have {\small $\wh{\ga}
< \f{ \ln \f{2}{\de} } { \varphi(\vep) } <  \f{ (1 + \vep) \ln \f{2}{\de} } { (2 \ln 2 - 1) \vep^2 } <  \f{ 4 (e
- 2) (1 + \vep) \ln \f{2}{\de} }{ \vep^2 }$}.

Second, we shall show $\wt{\ga} < \wh{\ga}$.  Clearly, if $\wt{\ga} < \f{1}{\vep}$, then $\wh{\ga} > \wt{\ga}$
is trivially true since $\wh{\ga} > \f{1}{\ga}$.  Thus, we can focus on the case that both $\wh{\ga}$ and
$\wt{\ga}$ are greater than $\f{1}{\vep}$.  For $\ga
> \f{1}{\vep}$, by Lemma \ref{lemA1}, we have $\eta < \ze$ and consequently $\wt{\mscr{Q}} (\vep, \ga) < \wh{\mscr{Q}} (\vep, \ga)$. As a result, $\de =
\wt{\mscr{Q}} (\vep, \wt{\ga}) = \wh{\mscr{Q}} (\vep, \wh{\ga})
> \wt{\mscr{Q}} (\vep, \wh{\ga})$.  Since $\wt{\mscr{Q}} (\vep, \wt{\ga}) > \wt{\mscr{Q}} (\vep, \wh{\ga})$ and
{\small $\f{ \pa [\wt{\mscr{Q}} (\vep, \ga) ] } { \pa \ga } < 0$}, we have $\wt{\ga} < \wh{\ga}$.

Third,  we shall show $(1 - \vep) \wh{\ga} < \wt{\ga}$.  In light of the fact that
\[
\eta = \f{ \vep \wh{\ga} - 1 } { \wh{\ga} - 1 }, \qqu \f{1}{1 - \vep} = \f{1}{1 - \ze} + \f{1}{\wh{\ga}},
\]
we have
\[
1 - \eta = (1 - \vep) \f{\wh{\ga}}{\wh{\ga} - 1}, \qu 1 - \ze = \f{ 1 - \vep } { 1 - \f{1 - \vep}{\wt{\ga}} },
\qu \f{ 1 - \eta  } { 1 - \ze  } = \f{\wh{\ga}}{\wh{\ga} - 1} \li ( 1 - \f{1 - \vep}{\wt{\ga}} \ri ).
\]
Therefore, if $(1 - \vep) \wh{\ga} \geq \wt{\ga}$, then $\f{ 1 - \eta  } { 1 - \ze  } \geq 1$, which implies
$\ze \leq \eta$. As a result, \bee \exp \li ( - \wt{\ga} \varphi(\vep)  \ri ) + \exp \li ( - \wt{\ga} \varphi(-
\ze) \ri ) & \geq & \exp \li ( - \wt{\ga} \varphi(\vep)
\ri ) + \exp \li ( - \wt{\ga} \varphi(- \eta)  \ri )\\
& > & \exp \li ( - \wh{\ga} \varphi(\vep)  \ri ) + \exp \li ( - \wh{\ga} \varphi(- \eta) \ri ) = \de, \eee which
contradicts  $\exp \li ( - \wt{\ga} \varphi(\vep)  \ri ) + \exp \li ( - \wt{\ga} \varphi(- \ze) \ri ) = \wt{
\mscr{Q} } (\vep, \wt{\ga}) = \de$.  Hence, it must be true that $(1 - \vep) \wh{\ga} < \wt{\ga}$.

Now, we shall show (\ref{tight0}). Since $\varphi(\vep) > 0$ and $\exp \li ( - \wt{\ga} \varphi(\vep) \ri ) <
\de = \wt{\mscr{Q}} (\vep, \wt{\ga})$,  we have $\f{\ln \f{1}{\de}}{ \varphi(\vep)} < \wt{\ga}$. Combining this
lower bound with the previously established upper bound yields \be \la{inq8} \f{\ln \f{1}{\de}}{ \varphi(\vep)}
< \wt{\ga} < \f{\ln \f{2}{\de}}{\varphi(\vep)}. \ee

Clearly,  as an immediate consequence of (\ref{inq8}), we have  {\small $\lim_{\de \to 0} \f{ \wt{\ga} } { \li [
\varphi(\vep) \ri ]^{-1} \ln \f{2}{\de} } = 1$}.   It remains to show {\small $\lim_{\vep \to 0} \f{ \wt{\ga}  }
{ \li [ \varphi(\vep) \ri ]^{-1} \ln \f{2}{\de} } = 1$}.  To this end, we need to prove {\small $\lim_{\vep \to
0} \f{ \exp ( - \wt{\ga} \varphi(\vep) ) \sh \de  } { \exp ( - \wt{\ga} \varphi(- \ze) ) \sh \de } = 1$}. It
suffices to show $\lim_{\vep \to 0} \wt{\ga} [\varphi(\vep) - \varphi(- \ze)] = 0$.  By virtue of (\ref{inq8})
and the condition $\f{1}{1 - \vep} = \f{1}{1 - \ze} + \f{1}{\wt{\ga}}$, \be \la{brid}
 \f{\ze}{\vep}  =  \f{1 -
\f{1 - \vep}{\vep \wt{\ga}} }{1 - \f{1 - \vep}{\wt{\ga}}} \geq \f{1 -  \f{1 - \vep}{\vep} \f{ \varphi(\vep) } {
\ln \f{1}{\de} }  }{1 - \f{1 - \vep}{\wt{\ga}}}. \ee Making use of the inequality of (\ref{brid}) and the facts
that  $\f{ \varphi(\vep) } { \vep } \to 0$ and $\wt{\ga} > \f{1 - \vep}{\vep} \to \iy$ as $\vep \to 0$, we have
\[
\liminf_{\vep \to 0} \f{\ze}{\vep} \geq \lim_{\vep \to 0} \f{1 -  \f{1 - \vep}{\vep} \f{ \varphi(\vep) } { \ln
\f{1}{\de} }  }{1 - \f{1 - \vep}{\wt{\ga}}} = 1.
\]
On the other hand, $\limsup_{\vep \to 0} \f{\ze}{\vep} \leq 1$ because $\ze < \vep$. Hence, \be \la{ineqB1}
\lim_{\vep \to 0} \f{\ze}{\vep} = 1. \ee Applying the upper bound in (\ref{inq8}),  we have
\[
\exp \li ( - \wt{\ga} \varphi(\vep)  \ri ) > \f{\de}{2} = \f{1}{2} \wt{ \mscr{Q} } (\vep, \wt{\ga})  = \f{1}{2}
\exp \li ( - \wt{\ga} \varphi(\vep)  \ri ) + \f{1}{2} \exp \li ( - \wt{\ga} \varphi(- \ze)  \ri ),
\]
which leads to $\exp \li ( - \wt{\ga} \varphi(\vep)  \ri ) >  \exp \li ( - \wt{\ga} \varphi(- \ze)  \ri )$, or
equivalently, \be \la{ineqB2} \varphi(\vep) - \varphi(- \ze) < 0. \ee By  (\ref{inq8}), (\ref{ineqB1}) and (
\ref{ineqB2} ), we have \bee
 \lim_{\vep \to 0} \f{\varphi(\vep) - \varphi(- \ze)}{ \varphi(\vep)}
& = & \lim_{\vep \to 0} \li [ 1 - \f{ \varphi(- \ze) }{ \ze^2 }
 \li ( \f{\ze}{\vep} \ri )^2 \f{\vep^2} { \varphi(\vep) } \ri ]\\
& = &  1 -  \lim_{\vep \to 0} \f{ \varphi(- \ze) }{ \ze^2 }
 \times \lim_{\vep \to 0} \li ( \f{\ze}{\vep} \ri )^2 \times \lim_{\vep \to 0} \f{\vep^2} { \varphi(\vep) }\\
& = &  1 - \f{1}{2} \times 1 \times 2  = 0 \eee and consequently, \bee
 \limsup_{\vep \to 0} \wt{\ga} [\varphi(\vep) - \varphi(- \ze)] & \leq & \limsup_{\vep
\to 0} \f{\ln \f{1}{\de}}{ \varphi(\vep)} [\varphi(\vep) - \varphi(- \ze)] = 0, \eee \bee
 \liminf_{\vep \to 0} \wt{\ga} [\varphi(\vep) - \varphi(- \ze)] & \geq & \liminf_{\vep
\to 0} \f{\ln \f{2}{\de}}{ \varphi(\vep)} [\varphi(\vep) - \varphi(- \ze)] = 0. \eee It follows that $\lim_{\vep
\to 0} \wt{\ga} [\varphi(\vep) - \varphi(- \ze)] = 0$ and thus {\small $\lim_{\vep \to 0} \f{ \exp ( - \wt{\ga}
\varphi(\vep) ) \sh \de  } { \exp ( - \wt{\ga} \varphi(- \ze) ) \sh \de } = 1$}.

Finally, since {\small $\f{ \exp ( - \wt{\ga} \varphi(\vep) )  } { \de } + \f{ \exp ( - \wt{\ga} \varphi(- \ze)
) } { \de } = 1$} and {\small $\lim_{\vep \to 0} \f{ \exp ( - \wt{\ga} \varphi(\vep) ) \sh \de  } { \exp ( -
\wt{\ga} \varphi(- \ze) ) \sh \de } = 1$},  we have {\small $\lim_{\vep \to 0} \f{ \exp ( - \wt{\ga}
\varphi(\vep) )  } { \de } = \f{1}{2}$}, which implies $\lim_{\vep \to 0} \f{ \wt{\ga} \varphi(\vep) } { \ln
\f{2}{\de} } = 1$.

\subsection{Proof of (VI) and (VII)}  First, we shall derive the upper bound of {\small $\Pr \li \{  \bs{n} \geq \f{\ga (1 + \vro) }{\mu}  \ri
\}$}. Since $\bs{n}$ is an integer, we have
\[
\Pr \li \{  \bs{n} \geq \f{\ga (1 + \vro) }{\mu}  \ri \}  = \Pr \li \{  \bs{n} \geq \li \lc \f{\ga (1 + \vro)
}{\mu} \ri \rc  \ri \} = \Pr \li \{  \bs{n} \geq \f{\ga (1 + \vro^*) }{\mu}  \ri \}
\]
where $\vro^*$ is a number satisfying {\small $\f{\ga (1 + \vro^*) }{\mu}  = \li \lc \f{\ga (1 + \vro) }{\mu}
\ri \rc$}. Clearly, $\vro^* \geq \vro$ by the definition of $\vro^*$. Let {\small $m = \f{\ga (1 + \vro^*)
}{\mu} - 1$}. Since $\vro
> \f{\mu}{\ga}$, we have $\f{\ga (1 + \vro) }{\mu} > 1$, which implies $m \geq 1$. Hence,
\[
\Pr \li \{  \bs{n} \geq \f{\ga (1 + \vro) }{\mu}  \ri \} = \Pr \{ \bs{n} \geq m + 1 \} = \Pr \{ X_1 + \cd + X_m
< \ga \} = \Pr \{ \ovl{X} < z \}
\]
with $\ovl{X} = \f{ \sum_{i=1}^m X_i } { m }$ and {\small $z = \f{\ga}{ m} = \f{\ga \mu}{ \ga (1 + \vro^*) -
\mu}$}. Note that $0 < z < \mu$ as a result of $\vro^* \ga \geq \vro \ga > \mu$.  It follows from Lemma
\ref{Hoe} that {\small \bee \Pr \li \{  \bs{n} \geq \f{\ga (1 + \vro) }{\mu} \ri \} & = & \Pr \{ \ovl{X} < z \}
\leq \exp (m z \mscr{M}( z, \mu))\\
& = & \exp \li ( \ga \mscr{M} \li ( \f{\ga \mu}{ \ga (1 + \vro^*) - \mu}, \mu \ri ) \ri ). \eee} Since $\vro^*
\geq \vro > \f{\mu}{\ga}$, applying Lemma \ref{decra}, we have $\mscr{M} \li ( \f{\ga \mu}{ \ga (1 + \vro^*) -
\mu}, \mu \ri ) \leq \mscr{M} \li ( \f{\ga \mu}{ \ga (1 + \vro) - \mu}, \mu \ri )$ and {\small \bee \Pr \li \{
\bs{n} \geq \f{\ga (1 + \vro) }{\mu} \ri \}
& \leq & \exp \li ( \ga \mscr{M} \li ( \f{\ga \mu}{ \ga (1 + \vro^*) - \mu}, \mu \ri ) \ri )\\
& \leq & \exp \li ( \ga \mscr{M} \li ( \f{\ga \mu}{ \ga (1 + \vro) - \mu}, \mu \ri ) \ri )\\
& = &  \exp \li (  \ga \li [ \ln \li ( 1 + \vro - \f{\mu}{\ga} \ri ) +
 \li (  \f{1 + \vro}{\mu} - \f{1}{\ga} - 1 \ri ) \ln \li ( \f{1 - \mu}{1 - \f{\ga \mu}
 { \ga (1 + \vro)  - \mu}} \ri ) \ri ] \ri )\\
& = &  \li ( 1 + \vro - \f{\mu}{\ga} \ri )^{ \li ( 1 + \vro - \f{\mu}{\ga} \ri ) \times \f{\ga}{\mu} }  \times
\li ( \f{1 - \mu} { 1 + \vro  - \f{\mu}{\ga} - \mu } \ri )^{ \li (1 + \vro - \f{\mu}{\ga} - \mu \ri ) \times
\f{\ga}{\mu} } \eee} for $\vro > \f{\mu}{\ga}$.

Now we bound {\small $\Pr \li \{ \bs{n}  \leq  \f{\ga (1 - \vro)}{\mu} \ri \}$}.  Invoking (\ref{remm}), we have
{\small \bee  \Pr \li \{ \bs{n} \leq \f{\ga}{\mu (1 + \vep)} \ri \}
 & \leq & \exp (\ga \mscr{M} (\mu (1 + \vep), \mu) )\\
& = & \exp \li ( \ga \li [ \ln \li ( \f{1}{1 + \vep} \ri ) + \li (\f{1}{\mu (1 + \vep)} - 1 \ri ) \ln \li ( \f{1
- \mu}{1 - \mu (1 + \vep) } \ri ) \ri ] \ri ) \eee} for $\mu (1 + \vep) < 1$.  Letting $\vro = \f{\vep}{1 +
\vep}$, we have  {\small \bee  \Pr \li \{ \bs{n}  \leq  \f{\ga (1 - \vro)}{\mu} \ri \}
 & \leq & \exp \li
( \ga \li [ \ln \li ( 1 - \vro \ri ) + \li (\f{1 - \vro}{\mu } - 1 \ri ) \ln \li ( \f{1 - \mu}{1 - \f{\mu} {1 -
\vro} } \ri ) \ri ]
\ri )\\
& = &   \exp \li ( \f{\ga}{\mu}  \li [ (1 - \vro) \ln \li ( 1 - \vro \ri ) + (1 - \vro - \mu) \ln \li ( \f{1 -
\mu}{1 - \vro  - \mu } \ri ) \ri ] \ri )\\
& = & \li ( 1 - \vro \ri )^{ (1 - \vro) \times  \f{\ga}{\mu} }  \times  \li ( \f{1 - \mu}{1 - \vro  - \mu } \ri
)^{ (1 - \vro - \mu)  \times  \f{\ga}{\mu}  } \eee} for $0 < \vro < 1 - \mu$.

\sect{Proof of Theorem 2}

We need the following preliminary result.

\beL \la{lem88} Let $\bs{k} = \bs{n} - \ga$.  Then, $\Pr \{ \bs{k} \geq s \} \leq \li ( \f{ s + \ga } { s } q
\ri )^s \li ( \f{ s + \ga } { \ga } p \ri )^\ga$ for $s > \bb{E} [ \bs{k}]$. Similarly, $\Pr \{ \bs{k} \leq s \}
\leq \li ( \f{ s + \ga } { s } q \ri )^s \li ( \f{ s + \ga } { \ga } p \ri )^\ga$ for $0 < s < \bb{E} [
\bs{k}]$. \eeL

\bpf

Note that $\bs{k}$ is a negative binomial random variable with distribution \[ \Pr \{ \bs{k} = k \} = \bi{\ga +
k - 1}{k} p^\ga q^k, \qqu k = 0, 1, 2, \cd.
\]
For any $t > 0$, {\small \bee \Pr \{ \bs{k} \geq s \} & = &
\Pr \{ e^{t ( \bs{k} - s )} \geq 1 \} \leq  \bb{E} \li [ e^{t ( \bs{k} - s )} \ri ] =  e^{- t s} \; \bb{E} \li [ e^{t \bs{k}} \ri ]\\
& = & e^{- t s} \; \sum_{k = 0}^\iy e^{t k } \bi{\ga + k -1}{k} p^\ga q^k\\
& = & e^{- t s} \; p^\ga \sum_{k = 0}^\iy  \bi{\ga + k -1}{k} (q e^t)^k  =  \phi(t) \eee} where $\phi(t) = e^{-
s} \li ( \f{p}{1 - q e^t} \ri )^\ga$. It can be checked that $\f{ d \ln \phi (t) }{d t} = - s + \f{ \ga qe^t} {
1 - qe^t } = 0$ if $s = (s + \ga)  qe^t$, in which $t = \ln \li( \f{s}{(s + \ga) q} \ri ) > 0$ because $s
> \bb{E} [ \bs{k}] = \f{q \ga}{p}$.  Substituting {\small $qe^t = \f{ s } { s + \ga}$} and {\small $e^{- t s} = \li (
\f{ s + \ga } { s } q \ri)^s$} in $\phi (t)$ yields the upper bound of $\Pr \{ \bs{k} \geq s \}$. Similarly, the
upper bound of $\Pr \{ \bs{k} \leq s \}$ can be established for $0 < s < \bb{E} [ \bs{k}]$.

\epf

Now we are in position to prove Theorem 2.  By the definition of the estimator $\wt{\bs{p}} = \f{\ga}{ \bs{k} +
\ga }$,
\[
\Pr \li \{  \li | \f{ \wt{\bs{p}} - p } { p } \ri | \geq \vep \ri \} =  \Pr \li \{   \bs{k} \leq \f{\ga}{p(1 +
\vep)} - \ga \ri \} + \Pr \li \{   \bs{k} \geq \f{\ga}{p(1 - \vep)} - \ga \ri \}.
\]
To bound {\small $\Pr \li \{   \bs{k} \leq \f{\ga}{p(1 + \vep)} - \ga \ri \} $}, we need to consider three cases
as follows.

\bed

\item (i): In the case of $p(1 + \vep)
> 1$, we have {\small $\Pr \li \{   \bs{k} \leq \f{\ga}{p(1 + \vep)} - \ga \ri \} = 0 < \exp \li ( - \ga
\varphi(\vep) \ri )$}.

\item (ii): In the case of  $p(1 + \vep) = 1$, we have {\small $\Pr \li \{   \bs{k} \leq \f{\ga}{p(1 + \vep)} - \ga \ri
\} = \Pr \{ \bs{k} = 0 \} = p^\ga  = (1 + \vep)^{- \ga} < \exp \li ( - \ga \varphi(\vep) \ri )$}.

\item (iii): In the case of $0 < p(1 + \vep) < 1$, applying Lemma \ref{lem88} with {\small $s = \f{\ga}{p(1 + \vep)} -
\ga < \f{\ga}{p} - \ga = \bb{E} [ \bs{k}]$}, we have {\small \bee \Pr \li \{   \bs{k} \leq \f{\ga}{p(1 + \vep)}
- \ga
\ri \} & \leq & \li ( \f{ s + \ga } { s } q \ri )^s \li ( \f{ s + \ga } { \ga } p \ri )^\ga\\
& = & \li ( \f{ \f{\ga}{p(1 + \vep)} } { \f{\ga}{p(1 + \vep)} - \ga } q \ri )^{\f{\ga}{p(1 + \vep)} - \ga }
\li ( \f{ \f{\ga}{p(1 + \vep)} } { \ga } p \ri )^\ga\\
& = & \li ( \f{ q } { 1 - p(1 + \vep) } \ri )^{\f{\ga}{p(1 + \vep)} - \ga }
\li ( \f{1}{1 + \vep}  \ri )^\ga\\
& = & \exp \li ( \ga \mscr{M}( (1 + \vep) p, p) \ri ). \eee} By the first statement of Lemma \ref{lembound}, we
have $\mscr{M}( (1 + \vep) p, p) \leq - \varphi( \vep )$.

\eed

Therefore, {\small $\Pr \li \{ \bs{k} \leq \f{\ga}{p(1 + \vep)} - \ga \ri \} \leq \exp \li ( - \ga \varphi( \vep
) \ri )$} is true for all cases.

To bound {\small $\Pr \li \{   \bs{k} \geq \f{\ga}{p(1 - \vep)} - \ga \ri \}$}, applying Lemma \ref{lem88} with
{\small $s = \f{\ga}{p(1 - \vep)} - \ga > \f{ \ga }{p} - \ga > \bb{E} [ \bs{k}] $} , we have {\small \bee
 \Pr \li \{
\bs{k} \geq \f{\ga}{p(1 - \vep)} - \ga \ri \} & \leq & \li ( \f{ s + \ga } { s } q \ri )^s \li ( \f{ s
+ \ga } { \ga } p \ri )^\ga\\
& = & \li ( \f{ q } { 1 - p(1 - \vep) } \ri )^{\f{\ga}{p(1 - \vep)} - \ga } \li ( \f{1}{1 - \vep}  \ri )^\ga\\
& = & \exp \li ( \ga \mscr{M}( (1 - \vep) p, p) \ri ). \eee} By the second statement of Lemma \ref{lembound}, we
have $\mscr{M}( (1 - \vep) p, p) \leq - \varphi( - \vep )$ and thus {\small $\Pr \{   \bs{k} \geq \f{\ga}{p(1 -
\vep)} - \ga \} \leq \exp \li ( - \ga \varphi(-\vep) \ri )$}.  Combining the two bounds of the tail distribution
probabilities of $\bs{k}$, we have
\[
\Pr \li \{  \li | \f{ \wt{\bs{p}} - p } { p } \ri | \geq \vep \ri \} \leq \exp \li ( - \ga \varphi(\vep) \ri ) +
\exp \li ( - \ga \varphi(-\vep)  \ri ) = \mscr{Q} (\vep, \ga).
\]
Clearly, $\mscr{Q} (\vep, \ga)$ is monotone decreasing with respect to $\ga$.  Because of such monotone property
and the fact that $\lim_{\ga \to \iy} \mscr{Q} (\vep, \ga) = 0, \; \lim_{\ga \to 0} \mscr{Q} (\vep, \ga)
> 1$, there exists a unique number $\ga^*$ such that  $\mscr{Q} (\vep, \ga^*) = \de$.

To derive the lower bound for $\ga^*$, applying Lemma \ref{lemcom}, we have $0 < \varphi(\vep) < \varphi(-\vep)$
and thus $2 \exp \li ( - \ga \varphi(- \vep) \ri ) < \mscr{Q} (\vep, \ga)$.  It follows that
\[
\max \{2 \exp \li ( - \ga^* \varphi(- \vep) \ri ), \; \exp \li ( - \ga^* \varphi(\vep) \ri ) \} < \de = \mscr{Q}
(\vep, \ga^*),
\]
from which we can obtain the lower bound of $\ga^*$.

To derive the upper bound for $\ga^*$, note that
\[
\mscr{Q} (\vep, \ga) - \wt{\mscr{Q}} (\vep, \ga)  = \exp \li ( - \ga \varphi(-\vep)  \ri ) - \exp \li ( - \ga
\varphi(-\ze)  \ri ) < 0 \]
 because $0 < \ze < \vep < 1$ and {\small $\f{d \varphi(- \vep) }{d \vep} > 0$} for
$0 < \vep < 1$.  Hence, $\de = \mscr{Q} (\vep, \ga^*) = \wt{\mscr{Q}} (\vep, \wt{\ga}) > \mscr{Q} (\vep,
\wt{\ga} )$. Since $\mscr{Q} (\vep, \ga)$ is monotone decreasing with respect to $\ga$, it must be true that
$\ga^* < \wt{\ga} < \f{ \ln \f{2}{\de} } { \varphi(\vep) }$.

Finally,  we consider the distribution of sample size.  Note that {\small
\[
\Pr \li \{ \bs{n} \geq \f{\ga (1 + \vro) }{p} \ri \} = \Pr \li \{ \bs{k} \geq s  \ri \} \] } where {\small $s =
(1 - p + \vro) \times \f{ \ga } {  p }$}. Note that {\small $\f{ s + \ga } { s } q = q \times \f{1 + \vro}{1 - p
+ \vro}$} and {\small $\f{ s + \ga } { \ga } p = 1 + \vro$}.  By Lemma \ref{lem88}, \bee \Pr \li \{ \bs{n} \geq
\f{\ga (1 + \vro) }{p} \ri \} & = & \Pr \li \{ \bs{k} \geq s \ri
\}\\
& \leq & \li ( \f{ s + \ga } { s } q \ri )^s \li ( \f{ s + \ga } { \ga } p \ri )^\ga\\
& = & \li (  q \times \f{1 + \vro}{1 - p + \vro} \ri
)^{ (1 - p + \vro) \times \f{ \ga } {  p } } \li (  1 + \vro \ri )^\ga \\
& = & \li (  \f{q }{1 - p + \vro} \ri )^{ (1 - p + \vro) \times \f{ \ga } {  p } } \li (  1 + \vro \ri )^{ \ga +
(1
- p + \vro) \times \f{ \ga } {  p } } \\
& = & \li (  \f{1 - p }{1 - p + \vro} \ri )^{ (1 - p + \vro) \times \f{ \ga } {  p } } \li (  1 + \vro \ri )^{
(1 + \vro) \times \f{ \ga } {  p } }. \eee Changing the sign of $\vro$ to negative yields
\[
\Pr \li \{ \bs{n} \leq \f{\ga (1 - \vro) }{p} \ri \}  \leq \li (  \f{1 - p }{1 - p - \vro} \ri )^{ (1 - p -
\vro) \times \f{ \ga } {  p } } \li (  1 - \vro \ri )^{ (1 - \vro) \times \f{ \ga } {  p } }.
\]
This completes the proof of Theorem 2.

\sect{Proof of Theorem 3}

For simplicity of notations, define {\small $C(p)  =  \Pr \li \{ \li | \wh{\bs{p}} - p \ri | < \vep  p \ri \}$}
and {\small $S(\ga, g, h,p) = \sum_{i = g}^h \bi{\ga + i - 1}{i} p^\ga (1- p)^i$}.  Then,
\bee C(p) & = & \Pr \li \{ \li | \f{ \ga - 1 }{\bs{k} + \ga - 1}- p \ri | < \vep p \ri \} \\
& = & \Pr \li \{ \f{\ga - 1 }{ (1 + \vep) p } - \ga + 1 <  \bs{k} < \f{\ga - 1  }{ (1 - \vep) p } - \ga + 1 \ri
\}\\
& = & \Pr \li \{ g(p) \leq \bs{k} \leq h(p) \ri \}  = S(\ga, g(p), h(p), p) \eee where
\[
g(p) = \li \lf \f{\ga - 1 }{ (1 + \vep) p } \ri \rf - \ga + 2, \qqu h(p) = \li \lc \f{\ga - 1  }{ (1 - \vep) p }
\ri \rc - \ga.
\]
It should be noted that $C(p), \; g(p)$ and $h(p)$ are actually multivariate functions of $p, \; \vep$ and
$\ga$.

We need some preliminary results.

\beL \la{minus} Let $p_\ell = \f{ \ga - 1 } { (1 - \vep) ( \ell + \ga - 1 ) }$ where $\ell \in \bb{Z}$. Then,
$h(p) = h(p_{\ell + 1}) = \ell$ for any $p \in (p_{\ell +1}, p_\ell)$. \eeL

\bpf Note that \[
 h(p)  =  \li \lc \f{\ga - 1  }{ (1 - \vep) p } \ri \rc - \ga  =  \li \lc \f{\ga - 1  }{ (1 - \vep) p_\ell  } \f{ p_\ell  } {p} \ri \rc - \ga\\
 =  \li \lc ( \ell + \ga - 1 ) \f{ p_\ell  } {p} \ri \rc - \ga.
 \]
Since {\small $1 < \f{ p_\ell  } {p} < \f{ p_\ell  } {p_{\ell + 1}} = \f{ \ell + \ga } { \ell + \ga - 1 }$} for
$p \in ( p_{\ell + 1}, p_\ell)$, we have
\[
\ell + \ga - 1 < \li \lc ( \ell + \ga - 1 ) \f{ p_\ell  } {p} \ri \rc \leq \li \lc ( \ell + \ga - 1 ) \f{ \ell +
\ga } { \ell + \ga - 1 } \ri \rc = \ell + \ga.
\]
Hence, $\ell - 1 < h(p) \leq \ell$.  Since $h(p)$ is an integer, it must be true that $h(p) = \ell =
 h(p_{\ell
+ 1})$.

\epf

\beL \la{plus} Let $p_\ell = \f{ \ga - 1 } { (1 + \vep) ( \ell + \ga - 1 ) } $ where $\ell \in \bb{Z}$. Then,
$g(p) = g(p_{\ell})=  \ell + 1$ for any $p \in (p_{\ell +1}, p_\ell)$. \eeL

 \bpf

 Note that \[
 g(p)  =  \li \lf \f{\ga - 1  }{ (1 + \vep) p } \ri \rf - \ga + 2 =
  \li \lf \f{\ga - 1  }{ (1 + \vep) p_\ell  } \f{ p_\ell  } {p} \ri \rf - \ga + 2\\
 =  \li \lf ( \ell + \ga - 1 ) \f{ p_\ell  } {p} \ri \rf - \ga + 2.
 \]
Since {\small $1 < \f{ p_\ell  } {p} < \f{ p_\ell  } {p_{\ell + 1}} = \f{ \ell + \ga } { \ell + \ga - 1 }$} for
$p \in ( p_{\ell + 1}, p_\ell)$, we have
\[
\ell + \ga - 1 \leq \li \lf ( \ell + \ga - 1 ) \f{ p_\ell  } {p} \ri \rf <  ( \ell + \ga - 1 ) \f{ \ell + \ga }
{ \ell + \ga - 1 }  = \ell + \ga.
\]
Hence, $\ell + 1 \leq g(p) < \ell + 2$.  Since $g(p)$ is an integer, it must be true that $g(p) = \ell + 1 =
 g(p_{\ell})$.

\epf

\beL \la{constant} Let $\al < \ba$ be two consecutive elements of the ascending arrangement of all distinct
elements of {\small $\{a, b \} \cup
 \{ \f{ \ga - 1 } { (1 - \vep) ( \ell + \ga - 1 ) }  \in (a, b) : \ell \in \bb{Z} \} \cup
 \{ \f{ \ga - 1 } { (1 + \vep) ( \ell + \ga - 1 ) } \in (a, b) : \ell \in \bb{Z} \} $}.
Then, both $g(p)$ and $h (p)$ are constants for any $p \in (\al, \ba)$.
 \eeL

 \bpf
Since $\al$ and $\ba$ are two consecutive elements of the ascending arrangement of all distinct elements of the
set, it must be true that there is no integer $\ell$ such that $\al < \f{ \ga - 1 } { (1 - \vep) ( \ell + \ga -
1 ) } < \ba$ or $\al < \f{ \ga - 1 } { (1 + \vep) ( \ell + \ga - 1 ) } < \ba$. It follows that there exist two
integers $\ell$ and $\ell^\prime$ such that {\small $(\al, \ba) \subseteq \li ( \f{ \ga - 1 } { (1 - \vep) (
\ell + \ga ) }, \f{ \ga - 1 } { (1 - \vep) ( \ell + \ga - 1 ) } \ri )$} and {\small $(\al, \ba) \subseteq \li (
\f{ \ga - 1 } { (1 - \vep) ( \ell^\prime + \ga ) }, \f{ \ga - 1 } { (1 - \vep) ( \ell^\prime + \ga - 1 ) } \ri
)$.} Applying Lemma \ref{minus} and Lemma \ref{plus}, we have {\small $g(p) = g \li ( \f{ \ga - 1 } { (1 - \vep)
( \ell + \ga - 1 ) } \ri )$} and {\small $h(p) = h \li ( \f{ \ga - 1 } { (1 - \vep) ( \ell^\prime + \ga ) } \ri
)$} for any $p \in (\al, \ba)$.

 \epf

\beL \la{lem_lim}
 For any $p \in (0,1)$, {\small $\lim_{t \downarrow 0} C(p + t) \geq C(p)$}
 and {\small $\lim_{t \downarrow 0} C(p - t) \geq C(p)$}.
\eeL

\bpf Observing that $g(p + t) \leq g(p)$ for any $t > 0$ and that
 \[ h(p + t)  =  h(p) + \li \lc \f{\ga - 1  }{ (1 - \vep) (p + t) }
- \li \lc \f{\ga - 1  }{ (1 - \vep) p } \ri \rc \ri \rc  = h(p) \] for {\small $- 1 < \f{\ga - 1  }{ (1 - \vep)
(p + t) }  - \li \lc \f{\ga - 1  }{ (1 - \vep) p } \ri \rc <
 0$},
i.e., {\small $0 < t < \f{\ga - 1  }{ (1 - \vep)  } \li ( \li \lc \f{\ga - 1  }{ (1 - \vep) p } \ri \rc - 1 \ri
)^{-1} - p$}, we have \be \la{ineqa} S(\ga, g(p + t), h (p + t), p + t ) \geq S(\ga, g(p), h (p), p + t ) \ee
for {\small $0 < t < \max \li \{ 1, \; \f{\ga - 1 }{ (1 - \vep)  } \li ( \li \lc \f{\ga - 1  }{ (1 - \vep) p }
\ri \rc  - 1 \ri )^{-1} \ri \} - p$}.  Since {\small \[ g(p + t)  =  g(p) + \li \lf \f{\ga - 1  }{ (1 + \vep) (p
+ t) } -
 \li \lf \f{\ga - 1  }{ (1 + \vep) p } \ri \rf \ri \rf, \]}  we have {\small \[ g(p + t) = \bec g(p) - 1 & \tx{for}
\; \f{\ga - 1  }{ (1 + \vep) p } =
 \li \lf \f{\ga - 1  }{ (1 + \vep) p } \ri \rf  \; \& \; 0 < t \leq \f{\ga - 1  }{ (1 + \vep)  }
 \li ( \li \lf \f{\ga - 1  }{ (1 + \vep) p } \ri \rf - 1
\ri )^{-1} - p,\\
g(p) & \tx{for} \; \f{\ga - 1  }{ (1 + \vep) p } \neq
 \li \lf \f{\ga - 1  }{ (1 + \vep) p } \ri \rf  \;
\& \; 0 < t \leq \f{\ga - 1  }{ (1 + \vep)  }
 \li ( \li \lf \f{\ga - 1  }{ (1 + \vep) p } \ri \rf
\ri )^{-1} - p. \eec
\]}
It follows that both $g(p + t)$ and $h(p + t)$ are independent of $t$ if $t > 0$ is small enough. Since $S(\ga,
g, h, p + t)$ is continuous with respect to $t$ for fixed $g$ and $h$, we have that $\lim_{t \downarrow 0}
S(\ga, g(p + t), h (p + t), p + t )$ exists.  As a result, \bee \lim_{t \downarrow 0} C(p + t) & = & \lim_{t
\downarrow 0} S(\ga,
g(p + t), h (p + t), p + t )\\
& \geq & \lim_{t \downarrow 0} S(\ga, g(p), h (p), p + t ) = S(\ga,  g(p), h (p), p ) = C(p), \eee where the
inequality follows from (\ref{ineqa}).

Observing that $h(p - t) \geq h(p)$ for any $t > 0$ and that
 \[ g(p - t)  =   g(p) + \li \lf \f{\ga - 1  }{ (1 + \vep) (p - t) } -
 \li \lf \f{\ga - 1  }{ (1 + \vep) p } \ri \rf \ri \rf  = g(p) \] for
 $0 < t < p - \f{\ga - 1  }{ (1 + \vep)  } \li ( 1 + \li \lf \f{\ga - 1  }{ (1 + \vep) p }
\ri \rf \ri )^{-1}$, we have \be \la{ineqb} S(\ga, g(p - t), h (p - t), p - t ) \geq S(\ga, g(p), h (p), p - t )
\ee for $0 < t < p - \f{\ga - 1  }{ (1 + \vep)  } \li ( 1 + \li \lf \f{\ga - 1 }{ (1 + \vep) p } \ri \rf \ri
)^{-1}$. Since
\[ h(p - t)  = h (p) + \li \lc \f{\ga - 1  }{ (1 - \vep) (p - t) } -
 \li \lc \f{\ga - 1  }{ (1 - \vep) p } \ri \rc \ri \rc,
\]
we have {\small \[ h(p - t) = \bec h(p) + 1  & \tx{for} \; \f{\ga - 1  }{ (1 - \vep)p } = \li \lc \f{\ga - 1  }{
(1 - \vep) p } \ri \rc \; \& \; 0 < t < p - \f{\ga - 1  }{ (1 - \vep)  } \li ( 1 + \li \lc
\f{\ga - 1}{ (1 - \vep) p } \ri \rc \ri )^{-1},\\
h(p) & \tx{for} \; \f{\ga - 1  }{ (1 - \vep)p } \neq \li \lc \f{\ga - 1 }{ (1 - \vep) p } \ri \rc \; \& \; 0 < t
< p - \f{\ga - 1  }{ (1 - \vep)  } \li ( \li \lc \f{\ga - 1}{ (1 - \vep) p } \ri \rc \ri )^{-1}. \eec
\]}
It follows that both $g(p - t)$ and $h(p - t)$ are independent of $t$ if $t > 0$ is small enough. Since $S(\ga,
g, h, p - t)$ is continuous with respect to $t$ for fixed $g$ and $h$, we have that $\lim_{t \downarrow 0}
S(\ga, g(p - t), h (p - t), p - t )$ exists. Hence, \bee \lim_{t \downarrow 0} C(p - t) & =
& \lim_{t \downarrow 0} S(\ga, g(p - t), h (p - t), p -t )\\
& \geq & \lim_{t \downarrow 0} S(\ga, g(p), h (p), p -t ) = S(\ga, g(p), h (p), p ) = C(p), \eee where the
inequality follows from (\ref{ineqb}).

\epf

\beL \la{unimodal} Let $0 < u < v < 1, \; h \geq 0$ and $g \leq h$. Then,
\[ \min_{p \in [u, v]} S(\ga, g, h, p) = \min \{ S(\ga, g, h, u), \;
S(\ga, g, h, v) \}.
\]
 \eeL

 \bpf

Since $\Pr \{ \bs{k} \leq k  \} = I_p(\ga, k + 1)$ where $I_p$ is the regularized incomplete beta function
\[
I_p (a, b) = \f{ B(p, a, b) } {  B(a, b) } = \f{ \int_0^p t^{a - 1} (1 - t)^{b - 1} d t  } { \int_0^1 t^{a - 1}
(1 - t)^{b - 1} d t   },
\]
we have {\small \be \la{ddd} \f{ \pa [\Pr \{ \bs{k} \leq l \} ] } { \pa p } = \f{ p^{\ga -1} (1 - p)^l } {
B(\ga, l + 1) }
> 0
\ee} for any integer $l \geq 0$.   To show the lemma, it suffices to consider $2$ cases as follows.

Case (i): $g \leq 0 \leq  h$. In this case, $S(\ga, g, h, p)
 = S(\ga, 0, h, p)$, which is increasing as a result of (\ref{ddd}).

Case (ii): $0 < g \leq h$. By (\ref{ddd}), for two integers $0 \leq k < l$, {\small \bee \f{ \pa [\Pr \{ k <
\bs{k} \leq l \} ] } { \pa p } & = &
\f{ p^{\ga -1} (1 - p)^l } { B(\ga, l + 1) } - \f{ p^{\ga -1} (1 - p)^k } { B(\ga, k + 1) }\\
&  = & \f{ p^{\ga -1} (1 - p)^k } { B(\ga, l + 1) } \li [ (1 - p)^{l - k} - \f{ l! (\ga + k)! } { k! (\ga + l)!
} \ri ]
> 0
\eee} if {\small $p <  1 - \li [ \f{ l! (\ga + k)! } { k! (\ga + l)! } \ri ]^{\f{1}{l - k}}$.  If follows that
$S(\ga, g, h, p)$} is a unimodal function of $p$.

From such investigation of the derivative of $C(p) = S(\ga, g, h, p)$ with respective to $p$, we can see that
one of the following three cases must be true: (1) $C(p)$ decreases monotonically for $p \in [u, v]$; (2) $C(p)$
increases monotonically for $p \in [u, v]$; (3) there exists a number $\se \in (u, v)$ such that $C(p)$
increases monotonically for $p \in [u, \se]$ and decreases monotonically for $p \in (\se, v]$.  It follows that
the lemma must be true for all cases.

 \epf

\beL \la{inbetween} Let $\al < \ba$ be two consecutive elements of the ascending arrangement of all distinct
elements of $\{a, b \} \cup
 \{ \f{ \ga - 1 } { (1
- \vep) ( \ell + \ga - 1 ) } \in (a, b) : \ell \in \bb{Z} \} \cup
 \{ \f{ \ga - 1 } { (1 + \vep) ( \ell + \ga - 1 ) } \in (a, b) : \ell \in \bb{Z} \} $.
Then,
 $C(p) \geq \min \{ C(\al), \; C(\ba) \}$ for any $p \in (\al, \ba)$.
\eeL

\bpf

By Lemma \ref{constant}, $g(p)$ and $h(p)$ are constants for any $p \in (\al, \ba)$. Hence, we can drop the
argument $p$ and write $g(p) = g, \; h(p) = h$ and $C(p) = S(\ga, g, h, p)$.

For $p \in (\al, \ba)$, define interval $[\al + t, \ba - t]$ with {\small $0 < t < \min \li (p - \al, \ba - p,
\f{\ba - \al}{2} \ri )$}. Then, $p \in [\al + t, \ba - t]$. By Lemma \ref{unimodal},
\[
C(p) \geq \min_{\mu \in [\al + t, \ba - t]} C(\mu) = \min \{ C(\al + t), \; C(\ba - t) \}
\]
for {\small $0 < t < \min \li (p - \al, \ba - p, \f{\ba - \al}{2} \ri )$}.  By Lemma \ref{lem_lim}, both
$\lim_{t \downarrow 0} C(\al + t)$ and $\lim_{t \downarrow 0}C(\ba - t)$ exist and are bounded from below by
$C(\al)$ and $C(\ba)$ respectively. Hence, {\small \bee C(p)  & \geq & \lim_{t \downarrow 0} \; \min \{ C(\al +
t), \; C(\ba - t) \} \\
& =  & \min \li \{ \lim_{t \downarrow 0} C(\al + t), \; \lim_{t \downarrow 0} C(\ba - t) \ri \} \geq \min \{
C(\al), \; C(\ba) \} \eee} for any $p \in (\al, \ba)$. \epf

\bsk

Finally, we can readily deduce Theorem 3.   The first statement on the minimum of the coverage probability
follows immediately from Lemma \ref{inbetween}.  The second statement on the minimum of the coverage probability
can be proved in a similar way.

\section{The Incomplete Work of Dagum et al.}

Let $Z_1, \; Z_2, \cdots$ be a sequence of i.i.d. random variables defined on the same probability space such
that $Z_i \in [0,1]$ and $\bb{E} [Z_i] = \mu_Z \in (0,1)$.
 In order to estimate $\mu_Z$, Dagum et al. proposed (in Section 2.1, page 1486 of \cite{Dagum})
 the following Stopping
Rule Algorithm:

\smallskip

{\small

Initialize $N \arl 0, \;\; S \arl 0$.

While $S <\Upsilon_1$ do: $N \arl N+1, \;\; S \arl S + Z_N$.

Return $\wh{\mu}_Z = \f{\Upsilon_1}{N}$ as the estimate of $\mu_Z$.}

\bigskip

In Section 2.1, page 1486 of \cite{Dagum}, Dagum et al. claimed that the reliability of the estimate
$\wh{\mu}_Z$ is asserted by the following ``Stopping Rule Theorem''.

\beT  Let $\vep \in (0,1)$ and $\de \in \li (0, 1 \ri )$. Then, $\Pr \{ | \wh{\mu}_Z - \mu_Z | \leq \vep \mu_Z
\} \geq 1 - \delta$.  \eeT

We would like to point out that the proof of Dagum et al. is not complete. There exists a significant gap which
cannot be patched by using their argument.

\subsection{The Proof of ``Stopping Rule Theorem'' by Dagum et al}

To exhibit the fallacy of the argument by Dagum et al., we shall represent their proof for ``Stopping Rule
Theorem'' in this section.   The following preliminary result is first established by Dagum et al. as Lemma 4.6
in page 1489 of \cite{Dagum}.

\begin{lemma} \la{lem4}
Let $\lambda = e - 2$. Let $\rho_Z = \max \{\sigma_Z^2, \; \vep \mu_Z  \}$ where $\sigma_Z^2$ is the variance of
$Z$.  Define $\xi_k = \sum_{i = 1}^k (Z_i - \mu_Z)$ for $k = 1, 2, \cd$.  Then, for any fixed $N
> 0$ and any $\beta \in [0, 2 \lambda \rho_Z]$, {\small \be \la{lem4ineq1} \Pr \li\{ \f{\xi_N}{N} \geq \beta  \ri\}
\leq \exp \li( \f{- N \beta^2} { 4 \lambda \rho_Z } \ri) \ee and \be \la{lem4ineq2}
 \Pr \li\{ \f{\xi_N}{N} \leq -\beta  \ri\} \leq \exp \li( \f{- N
\beta^2} { 4 \lambda \rho_Z } \ri). \ee}
\end{lemma}

\bigskip

The argument of Dagum et al. (in Section 5, page 1490-1491 of \cite{Dagum}) for the ``Stopping Rule Theorem''
proceeds as follows.  Let $N_Z$ be the sample size at the stopping time. Recall that $\wh{\mu}_Z =
\f{\Upsilon_1} {N_Z}$. It suffices to show that {\small \be \la{upper}
 \Pr \li \{N_Z <  \f{ \Upsilon_1 } { \mu_Z (1+ \vep) } \ri\} \leq
\f{\delta}{2} \ee} and that (equation (8),  page 1491 of \cite{Dagum}) {\small \be \la{lower}
 \Pr \li \{N_Z >  \f{ \Upsilon_1 } { \mu_Z (1- \vep) } \ri\} \leq
\f{\delta}{2}.  \ee} To show (\ref{upper}),  it suffices to consider the case that $\mu_Z (1+ \vep) \leq 1$,
since the theorem is trivially true if $\mu_Z (1+ \vep) > 1$. Let {\small $L = \li \lfloor \f { \Upsilon_1 } {
\mu_Z (1+ \vep) } \ri \rfloor$}. By the definitions of $\Upsilon_1$ and $L$,  {\small \be L  = \li \lfloor \f {
1 + (1+ \vep) \f{4 \lambda}{\vep^2} \;
 \ln \f{2}{\delta} } { \mu_Z (1+ \vep) } \ri \rfloor >  \f { 1 + (1+ \vep) \f{4 \lambda}{\vep^2} \;
 \ln \f{2}{\delta} } { \mu_Z (1+ \vep) } - 1  \geq \f{4 \lambda}{\vep^2 \mu_Z} \;
 \ln \f{2}{\delta} \label{bb}. \ee}
Since $N_Z$ is an integer, $N_Z < \f { \Upsilon_1 } { \mu_Z (1+ \vep) }$ implies $N_Z \leq L$. But $N_Z \leq L$
if and only if $S_L \geq \Upsilon_1$.  Thus, {\small $\Pr \li \{N_Z <  \f{ \Upsilon_1} { \mu_Z (1+ \vep) } \ri
\} \leq \Pr \{ N_Z \leq L \} = \Pr \{ S_L \geq \Upsilon_1 \}$} where $S_L = \sum_{i = 1}^L Z_i$.  Let $\beta =
\f{ \Upsilon_1} { L } - \mu_Z$.  Then, $\Pr \{ S_L \geq \Upsilon_1 \} = \Pr \{ S_L - \mu_Z L - \beta L \geq 0\}
= \Pr \li \{ \f{ \xi_L } { L } \geq \beta \ri \}$. Noting that $\vep \mu_Z \leq \beta \leq 2 \lambda \rho_Z$,
Lemma \ref{lem4} implies that {\small $\Pr \li \{ \f{ \xi_L } { L } \geq \beta \ri \} \leq \exp \li( \f{ - L
\beta^2 } {4 \lambda \rho_Z} \ri) \leq \exp \li( \f{ - L (\vep \mu_Z)^2 } {4 \lambda \rho_Z} \ri)$.}  Using the
last inequality of (\ref{bb}) and noting that $\rho_Z \leq \max \{ \mu_Z ( 1- \mu_Z), \; \vep \mu_Z  \} \leq
\mu_Z$, it follows that {\small $\Pr \li \{N_Z <  \f{ \Upsilon_1 } { \mu_Z (1+ \vep) } \ri\} \leq \Pr \li \{ \f{
\xi_L } { L } \geq \beta \ri \} \leq  \exp \li( \f{ - L (\vep \mu_Z)^2 } {4 \lambda \rho_Z} \ri) \leq \f{\de} {
2}$}. This completes the proof of (\ref{upper}).

Finally, instead of giving detailed argument in \cite{Dagum}, Dagum et al. claimed that the proof of
(\ref{lower}) is similar.

\subsection{A Hole in the Proof of Dagum et al}

We would like to point out that the proof of Dagum et al. is not complete because (\ref{lower}) cannot be shown
by a similar argument of Dagum et al. in proving (\ref{upper}). The gap is exhibited as follows.   To show
{\small $\Pr \li \{N_Z >  \f{ \Upsilon_1 } { \mu_Z (1- \vep) } \ri\} \leq \f{\delta}{2}$} in the similar spirit
of the first part, it is expected to construct an integer $L$ and a real number {\small $\ba  = \mu_Z - \f{
\Upsilon_1} { L }$} such that {\small \bel \la{c1} \li \{N_Z
> \f{ \Upsilon_1} { \mu_Z (1 - \vep) } \ri \} & \subseteq & \{ S_L \leq \Upsilon_1 \},\\
  \la{c3} \vep \mu_Z & \leq & \beta,\\
\la{c4} \ba & \leq & 2 \lambda \rho_Z,\\
\la{c2} L & \geq & \f{4 \lambda}{\vep^2 \mu_Z} \;
 \ln \f{2}{\delta} \eel} and consequently  {\small \bel
 \Pr \li \{N_Z >  \f{ \Upsilon_1} { \mu_Z (1 - \vep)
} \ri \} & \leq  & \Pr \{ S_L \leq \Upsilon_1 \} \la{c6}\\
&  = & \Pr \li \{ \f{ \xi_L } { L } \leq - \ba \ri \} \la{c7}\\
& \leq & \exp \li( \f{ - L \ba^2 } {4 \lambda \rho_Z} \ri) \la{c8}\\
& \leq & \exp \li( \f{ - L (\vep \mu_Z)^2 } {4 \lambda \rho_Z} \ri) \la{c9}\\
& \leq & \exp \li( \f{ - L \vep^2 \mu_Z } {4 \lambda} \ri) \la{c10}\\
& \leq & \f{\de} { 2} \la{c11}\eel} where (\ref{c6}) relies on (\ref{c1}); (\ref{c7}) is due to the definitions
of $\beta$ and $S_L$; (\ref{c8}) relies on (\ref{lem4ineq2}) of Lemma \ref{lem4} and (\ref{c4}); (\ref{c9})
relies on (\ref{c3}); (\ref{c10}) is due to the fact $\ro_Z \leq \max \{ \mu_Z (1 - \mu_Z), \; \vep \mu_Z \}
\leq \mu_Z$; (\ref{c11}) relies on (\ref{c2}).

Unfortunately, it is possible that (\ref{c3}) contradicts (\ref{c1})!  To see this, note that, by the definition
of $\ba$ and (\ref{c3}), we have $\vep \mu_Z \leq \mu_Z - \f{ \Up_1} { L }$,  i.e., {\small $L \geq \li \lc \f{
\Up_1 } { (1 - \vep) \mu_Z } \ri \rc$} since $L$ is an integer.  We can show that {\small $\li \{N_Z
> \f{ \Upsilon_1} { \mu_Z (1 - \vep) } \ri \} \subseteq \{ S_L \leq \Upsilon_1 \}$} is not true if {\small $L
\geq \li \lc \f{ \Up_1 } { (1 - \vep) \mu_Z } \ri \rc$}.  For this purpose, note that {\small $\li \{N_Z >  \f{
\Upsilon_1} { \mu_Z (1 - \vep) } \ri \} = \li \{N_Z > K \ri \} = \{ S_K < \Upsilon_1 \}$} where {\small $K = \li
\lf \f{ \Up_1 } { (1 - \vep) \mu_Z } \ri \rf$}.  For a random variable $Z$ with mean value $\mu_Z$ such that
$\f{ \Up_1 } { (1 - \vep) \mu_Z }$ is not an integer, we have {\small $K = \li \lf \f{ \Up_1 } { (1 - \vep)
\mu_Z } \ri \rf < \li \lc \f{ \Up_1 } { (1 - \vep) \mu_Z } \ri \rc \leq L$.} As a result, it is possible that
$\sum_{i=1}^K Z_i < \Upsilon_1 < \sum_{i=1}^L Z_i$, which implies that $\{ S_K < \Upsilon_1 \} \subseteq \{ S_L
\leq \Upsilon_1 \}$ is not true. Thus we have shown that (\ref{c1}) is not necessarily true if (\ref{c3}) is
satisfied. This demonstrates that  it is not possible to show {\small $\Pr \li \{N_Z
>  \f{ \Upsilon_1} { \mu_Z (1 - \vep) } \ri \} \leq \f{\de}{2}$} by using the argument of Dagum et al.

\sect{The Fallacy of Cheng's Reasoning}

In order to improve efficiency, Cheng revised the Stopping Rule Algorithm of Dagum et al. by replacing the
threshold value $\Up_1$ with a smaller number {\small $\bs{\al}$}.  See, pages 12--13, Algorithm 1 of Section 5,
and page 18, lines 9-10  of his paper \cite{Cheng}.

Cheng claimed that such a revised algorithm ensures $\Pr \{ \li | \wh{\mu}_Z - \mu_Z \ri | \leq \vep \mu_Z  \}
\geq 1 - \de$.  He first established Theorem 4 in page 7 of his paper, which is restated as Theorem \ref{CT}
follows.

\beT \la{CT} Let $Z_1, \cd, Z_n$ be i.i.d. random variables bounded in $[0,1]$ with common mean value $\mu_Z \in
(0,1)$. Let $0 < \vep < \min \{1, \; (1 - \mu_Z) \sh \mu_Z \}$. Then, {\small $\Pr \li \{ \li | \f{ \sum_{i=1}^n
Z_i }{n} - \mu_Z \ri | \leq \vep \mu_Z \ri \} \geq 1 - \de$} if {\small $n \geq \f{1}{\mu_Z} \f{1}{ (1 + \vep)
\ln (1 + \vep) - \vep } \ln \li ( \f{2}{\de} \ri )$}.  \eeT

In the first paragraph of page 18 of his paper \cite{Cheng}, after defining events {\small \[ E_1 = \li \{  0 <
\mu_Z < \f{ \wh{\mu}_Z } {1 + 3 \vep } \ri \}, \qu E_2 = \li \{ \f{ \wh{\mu}_Z } {1 + 3 \vep } \leq \mu_Z < \f{
\wh{\mu}_Z } {1 + 2 \vep } \ri \}, \]}
 {\small \[
 E_3 = \li \{ \f{ \wh{\mu}_Z } {1 + 2 \vep } \leq \mu_Z < \f{ \wh{\mu}_Z } {1 +  \vep
} \ri \}, \qu E_4 = \li \{ \mu_Z \geq \f{ \wh{\mu}_Z } {1 +  \vep } \ri \}, \]} Cheng applied the law of total
probability to write
\[
\Pr \li \{ \li | \wh{\mu}_Z - \mu_Z  \ri | \leq \vep \mu_Z \ri \} = \sum_{i = 1}^4 \Pr \{ \li | \wh{\mu}_Z -
\mu_Z  \ri | \leq \vep \mu_Z  \mid E_i \} \Pr \{ E_i \} \] and attempted to show that the right-hand side of the
equality is bounded from below by $1 - \de$.  Unfortunately, Cheng made a fundamental mistake in bounding $\Pr
\li \{ \li | \wh{\mu}_Z - \mu_Z \ri | \leq \vep \mu_Z \mid E_4  \ri \}$ and other terms alike.  He noted that
$S_{N_Z} = \sum_{i = 1}^{N_Z} X_i \geq \bs{\al}$ and thus {\small $N_Z \geq \f{1}{ \wh{\mu}_Z } \f{1}{ \ln (1 +
\vep) - \vep \sh (1 + \vep) }  \ln \li ( \f{2}{\de_s} \ri )$} because of the definitions of the stopping rule
and $\bs{\al}$. Conditioning upon $\mu_Z \geq \f{ \wh{\mu}_Z } {1 +  k \vep }$ with some constant $k$, he
obtained \be \la{upon} N_Z \geq \f{1}{ \mu_Z } \f{1 + \vep}{1 + k \vep} \f{\ln \li ( \f{2}{\de_s} \ri )}{ (1 +
\vep) \ln (1 + \vep) - \vep }  = \f{1}{ \mu_Z } \f{1}{ (1 + \vep) \ln (1 + \vep) - \vep } \ln \li ( \f{2}{\de_s}
\ri )^{\f{1 + \vep}{1 + k \vep}} \ee and then applied Theorem \ref{CT} to claim {\small \be \la{err} \Pr \li \{
\li | \f{ \wh{\mu}_Z - \mu_Z } { \mu_Z } \ri | \leq \vep  \mid \mu_Z \geq \f{ \wh{\mu}_Z } {1 + k \vep } \ri \}
\geq 1 - 2 \li (\f{\de_s}{2} \ri )^{ \f{1 + \vep}{1 + k \vep} }, \ee} which was equation (28) of page 17 in his
paper \cite{Cheng}.

Here Cheng made a subtle and critical mistake by {\it illegally} applying Theorem \ref{CT}.  The reason is that
{\it the sample size requirement of Theorem \ref{CT} is independent of samples, while the validity of
(\ref{upon}) depends on samples}. Consequently, (\ref{err}) is not justified.  This affects subsequent relevant
development.

\end{document}